\documentclass[a4paper,11pt,oneside]{article}

\usepackage[top=2cm, bottom=2cm, left=2cm, right=2cm]{geometry}

\usepackage{authblk}

\usepackage{amsmath,amsthm,amsfonts,amssymb,amscd}
\usepackage{subfig}
\usepackage{booktabs}
\usepackage{mathtools}
\usepackage{graphicx,bm,xcolor}
\usepackage{algorithm}
\usepackage{algpseudocode}
 \usepackage{url}
\usepackage{siunitx}
\usepackage{amsmath}
\usepackage{amsthm}
\usepackage{amsbsy}
\usepackage{amssymb}
\usepackage{enumerate}

\usepackage{tikz}
\usetikzlibrary{hobby}
\tikzset{every path/.append style={line width=1pt}}
\usepackage{latexsym}
\usepackage{accents}
\usepackage{amsfonts}
\usepackage{amstext}
\usepackage{amsxtra}
\usepackage{hyperref}
\usepackage{color}
\usepackage{graphicx}
\usepackage{mathtools}
\usepackage{mathrsfs}
\usepackage{autobreak}
\usepackage{cancel}
\usepackage{pgfplots}
\usepackage{pgfkeys}
\usepackage{tabu}   
\usepackage{comment}
\usepackage{fancyhdr}
\usepackage{algorithm}
\usepackage{algpseudocode}
\usepackage{nameref}
\usepackage{chngcntr}

    \def\addlegendimage{\csname pgfplots@addlegendimage\endcsname}
\pgfplotsset{
cycle list={%
{draw=black,densely dotted},
{draw=blue,densely dotted},
{draw=red,densely dotted},
{draw=green,densely dotted}}}
\usepgfplotslibrary{fillbetween}

\newtheorem{Satz}{Satz}[section]

\newtheorem{Lemma}[Satz]{Lemma}

\newtheorem{Remark}[Satz]{Remark}
\newtheorem{Problem}[Satz]{Problem}

\normalsize
\newcommand{\Rr}{\mathbb{R}}

\newcommand{\eps}{\varepsilon}

\newcommand{\integral}[4]{\int_{#1}^{#2} #3  \def\temp{#4}\ifx\temp\empty
    \empty
\else
  \,\mathrm{d}#4\,
\fi}

\newcommand{\uproman}[1]{\uppercase\expandafter{\romannumeral#1}}

\newcommand\restr[2]{{
  \left.\kern-\nulldelimiterspace 
  #1 
  \vphantom{\big|} 
  \right|_{#2} 
  }}

\pgfplotsset{compat=1.15}

\definecolor{green}{RGB}{60,179,113}

\newcommand\RA{\textcolor{black}}
\newcommand\RC{\textcolor{black}}
\newcommand\auth{\textcolor{black}}

\begin{document}

%
%
\title{A modified combined active-set Newton method for solving 
phase-field fracture into the monolithic limit}

\author[1]{Leon Kolditz}
\author[2]{Katrin Mang}
\author[1,3]{Thomas Wick}

\affil[1]{Leibniz Universität Hannover, Institute of Applied Mathematics, Welfengarten 1, 30167 Hannover, Germany}
\affil[2]{Leibniz Universität Hannover, Institute of Structural Analysis, Appelstrasse 9A, 30167 Hannover, Germany}
\affil[3]{Universit\'e Paris-Saclay, LMPS - Laboratoire de Mecanique Paris-Saclay, 91190 Gif-sur-Yvette, France}

\date{}

\maketitle

\begin{abstract}
In this work, we examine a numerical phase-field fracture framework 
in which the crack irreversibility constraint is treated with a primal-dual active set method and a linearization is used in the degradation function to enhance the numerical stability.
The first goal is to carefully derive from a complementarity system our primal-dual active set formulation, which 
has been used in the literature in numerous studies, but for phase-field fracture 
without its detailed mathematical derivation yet. Based on the latter, we formulate 
a modified combined active-set Newton approach that significantly reduces the computational cost in comparison to comparable prior algorithms for quasi-monolithic settings. 
For many practical problems, Newton converges fast, but active set needs many iterations, for which three different efficiency improvements are suggested in this paper.
Afterwards, we design an iteration on the linearization  
in order to iterate the problem
to the monolithic limit. Our new algorithms are implemented 
in the programming framework pfm-cracks [T. Heister, T. Wick; 
pfm-cracks: A parallel-adaptive framework for phase-field fracture propagation,
Software Impacts, Vol. 6 (2020), 100045].
In the numerical examples, we conduct performance studies and investigate efficiency enhancements. 
The main emphasis is on the cost complexity by keeping the accuracy of numerical solutions and goal functionals. Our algorithmic suggestions are substantiated with the help of several benchmarks in two and three spatial dimensions. Therein, predictor-corrector adaptivity and parallel performance 
studies are explored as well.
\textbf{This manuscript has been accepted for publication in the Journal 
Computer Methods in Applied Mechanics and Engineering 
under the DOI \url{https://doi.org/10.1016/j.cma.2023.116170}.
}
\end{abstract}
	



\section{Introduction}\label{sec:introduction}
This work is devoted to the cost-efficient numerical quasi-monolithic solution 
of phase-field fracture problems. Since the pioneering studies from
\cite{francfort1998revisiting,bourdin2000revisited} and \cite{MieWelHof10a,KuMue10}, 
the variational phase-field approach to fracture has gained a lot of attention.
Summaries and overview monographs include 
\cite{BourFraMar08,AmGeraLoren15,WuNgNgSuBoSi19,wick2020multiphysics,BouFra19,DiLiWiTy22},
\RC{and recent current trends are outlined in selected chapters in \cite{AlHuSoWeWeiMa22}.}

Specifically, various algorithms and numerical 
studies on the nonlinear and linear solution have been undertaken to date.
In principle, iterative coupling (operator splitting; alternating minimization; \auth{staggered; partitioned}) and monolithic 
approaches are distinguished. In the early work, only iterative algorithms 
were used \cite{bourdin2000revisited}, for which convergence results could 
be established theoretically \cite{Bour07,BuOrSue10}, \RC{and extended to adaptive discretizations \cite{BuOrSue13}}. Furthermore, 
\cite{miehe2010phase} introduced as well operator splitting together 
with a strain history function for crack irreversibility. 
\RC{Recent promising
advancements in partitioned schemes include \cite{STORVIK2021113822,MesBouKhon15,LUO2023115787} and 
some of the following references do also consider partitioned approaches as we will mention.
}

The first monolithic solution methods date back to 
\cite{GeLo16,Wi17_CMAME,Wi17_SISC}. Concerning the 
arising linear systems within the nonlinear iterations, the first significant
study on linear solvers for phase-field fracture 
problems was done in \cite{FaMau17}.
More recently, further 
monolithic approaches were implemented in \cite{KoKr20,KoKoKr23,LAMPRON2021114091, lampron2023phase,wambacq2021interior,KRISTENSEN2020102446}
and a truncated non-smooth Newton multigrid solver was 
proposed in \cite{GrKieSa20}. 
\RC{In \cite{MaViBo16,SINGH201614,BHARALI2022114927} an arc-length approach was introduced 
in order to enhance the robustness of the solution process. Specifically, in 
\cite{SINGH201614,wambacq2021interior} comparisons of 
monolithic and partitioned 
methods were undertaken. Moreover, in Abaqus \cite{ma14081913}, implementations with comparisons between monolithic and partitioned 
schemes have been considered over the last years \cite{NAVIDTEHRANI2021100050,LIU201635}.}
As investigated in \cite{wu2020bfgs}, classical Newton methods fail to converge for large incremental steps. Thus, monolithic solution approaches often extend and modify the Newton method to improve its stability.
In order to increase the robustness, quasi-monolithic approaches 
were first proposed in \cite{heister2015primal} in which in one single 
term, actually in the quasi-linear displacement equation, an explicit 
linearization of the phase-field variable is introduced, but otherwise the 
system is still solved in a monolithic fashion. 
\RA{This linearization is performed in two different ways: As a first approach, the solution of the previous incremental step replaces the unknown phase-field in the displacement equation. In a second approach, to enhance the accuracy, an extrapolation of two previous incremental steps is used as an approximation.}
This 
system is numerically very robust and a block-diagonal preconditioner 
with algebraic multigrid preconditioning yields an efficient 
parallel solution \cite{HeiWi18_pamm} with an open-access 
implementation in \cite{heister2020pfm-cracks}. This quasi-monolithic 
approach was also used in \cite{JoLaWi20} in which a geometric 
matrix-free multigrid method was employed for the arising linear systems.
\RC{Finally, for large-scale problems towards practical applications, multiscale 
or (non-intrusive) global-local approaches become important to keep the computational 
cost reasonable \cite{GeNoiiAllLo18}. Therein, local problems contain the high-fidelity full nonlinear problem, which 
is coupled to a simpler (often linear) global problem in the surrounding medium. Later, this 
idea was further extended to adaptive schemes \cite{NoAlWiWr19}, multilevel methods \cite{ALDAKHEEL2021114175}, extended/generalized finite elements for fracture \cite{GePlTuDo20},
and thin-walled large deformations with phase-field \cite{LIU2022115410}.
More recently, due to the success of physics-informed neural neural networks, phase-field 
fracture numerical solutions are also solved with those methods \cite{GOSWAMI2020102447,GOSWAMI2022114587}.
}

The starting point of the current work is the quasi-monolithic solution
proposed in \cite{heister2015primal}. Therein a combined 
Newton method was \auth{designed} in which the nonlinearities of 
constitutive laws and nonlinear coupling were combined with a primal-dual 
active set \RA{(PDAS)} iteration to satisfy the crack irreversibility condition.
First, \auth{in the current work} the mathematical derivation of the primal-dual active 
set method is investigated in much more detail than in \cite{heister2015primal},
by starting from a weak complementarity formulation.
This provides (for the first time) a rigorous mathematical justification in terms 
of phase-field fracture formulations. 
The second objective is a significant improvement of the efficiency of the 
combined Newton method since in numerous simulations it was observed that 
Newton's method converges fast, but the active set method needs many iterations, but 
often only a few degrees of freedom change. Sometimes, only the same degrees of 
freedom switch from active to non-active, which is known as cycling 
in the numerical optimization literature \cite{Curtis2015}. However, the physics 
of the solution does not change anymore significantly. To this end, 
we propose four test cases to \auth{find} the impact of the active set method. Within these investigations, we analyze the role of a constant parameter inside the active set method. This constant was already mentioned in \cite{hueber2005primaldual,popp2009finite,schroeder2016semismooth} and some theoretical statements were made, but it was not analyzed further.
\RC{
The third aim of this work is a bit specific. In \cite{heister2015primal}, we introduced 
a linearization in terms of some extrapolation of (e.g., two) previous time-step solutions of the phase-field variable in the displacement equation in order to 
obtain a block-triangular Newton system matrix. This procedure results into an extremely 
robust nonlinear and linear numerical solution in which the linear system can be preconditioned 
with algebraic multigrid \cite{HeiWi18_pamm,heister2020pfm-cracks} or geometric multigrid methods \cite{JoLaWi20}.
However, introducing such a linearization results into a temporal discretization error.
For stationary or slowly growing fractures, this discretization error is sufficiently small. 
However, the error becomes significant the faster the fracture grows; see e.g., \cite{Wi17_SISC}.
In order to reduce this temporal discretization error, one approach, proposed 
in \cite{wick2020multiphysics}, and further tested in \cite{LAMPRON2021114091}, 
is to apply an additional iteration on the current time-step solution. Thus, we shall 
further investigate these additional iterations in terms of computational cost versus 
accuracy, as well as its behavior when combined with adaptive mesh refinement.
}

Gathering these developments together with the existing 
features in \texttt{pfm-cracks}, namely MPI (message passing interface)
parallelization and predictor-corrector adaptive mesh refinement, 
a new final algorithm 
is designed by us in which the modified Newton method and the iteration on 
linearization (ItL) are realized and subsequently implemented. Various settings 
for stationary, nonstationary and two- and three-dimensional configurations 
are adopted to study the performance of our new algorithmic framework. 
In the end, our aim is to substantiate a significant advancement of 
a parallel-adaptive phase-field fracture framework with a modern 
primal-dual active set strategy accounting for inequality constraints 
and iterating into the monolithic limit.

The outline of this paper is as follows: In Section \ref{sec:problem_formulation}, we describe the underlying setting of phase-field fracture. Then, we 
derive a complementarity system and obtain a series of intermediate mathematical 
results. In Section \ref{sec:numerical_solution}, our modified combined 
Newton scheme is introduced with four different cases. Then, in Section \ref{sec:iteration_on_extra} an iteration is introduced, which iterates the 
linearized phase-field value into the monolithic limit. Afterwards, in Section \ref{sec:numerical_tests} several numerical tests are conducted in order to 
study our novel algorithms for stationary, nonstationary two and three dimensional test cases.
Our work is summarized in Section \ref{sec_conclusions}.

\section{Problem formulation}\label{sec:problem_formulation}

\subsection{Notations}\label{subsec:notations}
This section introduces basic notations and the coupled variational inequality system (CVIS) and its discretized formulation, 
which is later used for simulations. The scalar-valued $L^2$-product on a sufficiently smooth, 
bounded domain $G\subset \Rr^d$, $d=2,3$ is denoted by
\begin{equation*}
    (x,y)_{L^2(G)} \coloneqq \integral{G}{}{x\cdot y}{G},
\end{equation*}
whereas the vector-valued $L^2$-product is defined by
\begin{equation*}
    (X,Y)_{L^2(G)}  \coloneqq \integral{G}{}{X:Y}{G},
\end{equation*}
with the Frobenius product $X : Y$ of two vectors $X,Y$ and the material $\Omega \subset \Rr^d$, which is sufficiently smooth and bounded. 
If there is no subscript provided, the $L^2$-product over the whole domain $\Omega$ is meant.

The Euler-Lagrange equations, arising from directional derivatives of an energy functional \cite{francfort1998revisiting,bourdin2000revisited} (pure elasticity), i.e., 
\cite{MiWheWi13a,MiWheWi19} (pressured fractures in a monolithic setting), regularized with an Ambrosio-Tortorelli 
approximation \cite{ambrosio1990approximation, ambrosio1992freediscontinuity}, consist of a displacement equation and the phase-field inequality and arise in a weak formulation. 
Consequently as solution variables, we have the displacement function $u:\Omega \rightarrow \Rr^d$ and the phase-field function $\varphi : \Omega \rightarrow [0,1]$. \RA{It is defined such that $\varphi=1$ in the intact part of the domain, $\varphi = 0$ in the fully broken part of the domain and $0<\varphi<1$ in the transition zone.}
The continuous-level solution sets are defined as 
\begin{align*}
    \mathcal{V} & \coloneqq H^1_0(\Omega), \quad
    \mathcal{K}^n \coloneqq \{\psi \in \mathcal{W}\, \vert\, \psi - \varphi^{n-1} \leq 0\quad \text{a.e.\ in}\ \Omega\}, \quad\text{where }  \mathcal{W} \coloneqq H^1(\Omega).
\end{align*}
The convex set for the phase-field variable arises due to the crack irreversibility constraint $\partial_t \varphi \leq 0$, which 
reads in inremental form $\varphi^n \leq \varphi^{n-1}$, where $\varphi^n:=\varphi(t_n)$ and $\varphi^{n-1}:=\varphi(t_{n-1})$.
Here, $t_n$ goes from $n=0$ (initial condition) until $n=N$ (end time condition)
\RA{such that we have the incremental grid $t_0,\dots,t_N$ and the step size $k_n = t_n - t_{n-1}$. 
However, $k_n$ does only appear implicitly in the inequality constraint, and no} time derivatives arise 
in the two governing problem statements, \auth{therefore,} the resulting framework is of quasi-static fashion.

\subsection{The phase-field fracture Euler-Lagrange equations}
In this short section, we formulate a monolithic CVIS of phase-field fracture.
First, we introduce constitutive laws and material parameters.
We work with the classical stress tensor of linearized elasticity defined by
\begin{equation*}
    \sigma(u) = 2\mu e(u) + \lambda \operatorname{tr}(e(u))I,
\end{equation*}
with the Lam\'e parameters $\mu>0$, and $\lambda$ with $3\lambda + 2\mu >0$,
and the identity matrix $I$. The symmetric strain tensor $e(u)$ is given by
\begin{equation*}
    e(u)\coloneqq \frac{1}{2}\left( \nabla u + \nabla u^T \right).
\end{equation*}
Moreover, the critical energy release rate is denoted by $G_C$ with $G_C>0$.

Then, the Euler-Lagrange equations are given by 
\cite{MiWheWi13a,MiWheWi19,wick2020multiphysics}
\begin{Problem}\label{problem:euler_lagrange_nonlinear}(Euler-Lagrange equations)
For some given initial value $\varphi^0$ and for the incremental steps $t_n$ with $n=1,...,N$, find $(u^n,\varphi^n) \in \mathcal{V} \times \mathcal{K}^n$ such that
\begin{align*}
    \left(g(\varphi^{\RA{n}})\sigma(u^n),e(\psi^u)\right) + ((\varphi^n)^2p^n, \operatorname{div} \psi^u) &= 0 \quad  \forall \psi^u \in \mathcal{V},
    \intertext{and}
    (1-\kappa)(\varphi^n\sigma(u^n) : e(u^n), \psi^\varphi - \varphi^n) +2(\varphi^n p^n \operatorname{div} u^n, \psi^\varphi-\varphi^n) \qquad &\quad \\
    \qquad + G_C\left(\frac{1}{\eps}(1-\varphi^n,\psi^\varphi-\varphi^n) + \eps(\nabla\varphi^n,\nabla (\psi^\varphi-\varphi^n))\right) &\geq 0 \quad \forall \psi^\varphi \in \mathcal{K}^n \cap L^\infty(\Omega),
\end{align*}
for a given pressure $p^n\in L^\infty(\Omega)$ (Sneddon's test for 
pressured fractures \cite{MiWheWi13a,MiWheWi19}) and \auth{$p^n \equiv 0$} for fracture in pure elasticity.
\end{Problem}
Therein, the degradation function is given by $g(\varphi^n) := (1-\kappa)(\varphi^n)^2+\kappa$. The bulk regularization parameter $\kappa$ is necessary to avoid irregularities in the system matrix. If the phase-field function $\varphi$ is $0$, we obtain zero-entries on the diagonal of the system matrix. To avoid this, we employ $\kappa>0$. We have to ensure that $\kappa$ is not too large since it yields \auth{a} perturbation of the \auth{physics of the} system, but it needs to be large enough to prevent irregularities. The second regularization parameter $\varepsilon$ appears due to the Ambrosio Tortorelli approximation \auth{\cite{ambrosio1990approximation,ambrosio1992freediscontinuity}}. The $\Gamma$-convergence theory, e.g \cite{Braides1998}, states on an energy-level that the regularized terms (under certain assumptions) converge to the underlying unregularized model as $\varepsilon \rightarrow 0$. A rigorous proof of pressurized phase-field fracture of the one dimensional case is done in \cite{EngSchu2016pressfracture}, a proof in higher dimensions is \auth{established} in \cite{sommer2019unfitted}.
For the discretized problem, we also have to require $h=o(\varepsilon)$ for the discretization parameter $h$. Finding an optimal setting for the regularization parameters $\varepsilon$ and $\kappa$, which yields the best compromise between computational cost and 
$\Gamma$-convergence theory needs some work and is highly test dependent \cite{kolditz2022relation}. 
In practice, \auth{we usually must choose} $h<\varepsilon$  
and \auth{often} realized with $\varepsilon = 2h$.

\subsection{A phase-field fracture formulation with degradation function linearization}\label{sec_PFF_linearization}
Problem \ref{problem:euler_lagrange_nonlinear} has several nonlinearities, namely 
nonlinear coupling of variables in the displacement PDE \auth{(partial differential equation)} and the phase-field 
inequality, while the inequality constraint introduces its own nonlinear behavior.
A brief analysis of the coupling terms reveals that the nonlinear 
behavior in the displacement equation is more severe (being of 
quasi-linear type; \auth{for the definition of quasi-linear, we refer the reader to} \cite{Evans2010}) in comparison to the phase-field 
part, which is semi-linear only. Therefore, if linearizations are of 
interest, it is reasonble to address the quasi-linear part first.
Indeed, the fully monolithic system is a big challenge to be solved 
\cite{GeLo16,Wi17_CMAME,Wi17_SISC,KoKr20,KoKoKr23,LAMPRON2021114091,wambacq2021interior,GrKieSa20,wu2020bfgs}. We follow our prior work \cite{heister2015primal} and formulate 
a quasi-monolithic system in which we linearize $(\varphi^n)^2$ by using known information about older incremental steps. In the first approach, we use an extrapolation for $\varphi^n$ such that
\begin{equation*}
    \Tilde{\varphi}^n \coloneqq  \Tilde{\varphi}^n(\varphi^{n-1},\varphi^{n-2}) = \varphi^{n-2} \frac{t_n-t_{n-1}}{t_{n-2}-t_{n-1}} + \varphi^{n-1} \frac{t_n-t_{n-2}}{t_{n-1}-t_{n-2}}.
\end{equation*}
In a second approach, we simply use the solution from the previous timestep such that 
\begin{equation*}
    \Tilde{\varphi}^n \coloneqq  \Tilde{\varphi}^n(\varphi^{n-1}) = \varphi^{n-1}.
\end{equation*}
Furthermore, the stress can be split into a compressive and a tensile part such that the energy degradation only acts on the tensile stress, \auth{which introduces a third nonlinearity in the system.} With these modifications, 
we obtain the following modified form of Problem \ref{problem:euler_lagrange_nonlinear}.
\begin{Problem}\label{problem:euler_lagrange_linearized}(Linearized Euler-Lagrange equations with stress-splitting)
For a given $\varphi^0$ and for every incremental step $t_n$ with $n=1,...,N$, find $U^n \coloneqq\{u^n,\varphi^n\} \in \mathcal{V} \times \mathcal{K}^{\auth{n}}$ such that it holds for $\Phi \coloneqq \{0,\varphi^n\} \in \mathcal{V}\times \mathcal{K}^n$
\begin{align*}
    A(U^n)(\Psi-\Phi) \geq 0 \quad \forall \Psi \coloneqq \{\psi^u,\psi^\varphi\} \in \mathcal{V}\times \mathcal{K}^n\cap L^\infty,
\end{align*}
where $A(U^n)(\Psi-\Phi)$ is defined as 
\begin{align*}
    A(U^n)(\Psi-\Phi) &\coloneqq\left(g(\Tilde{\varphi}^n)\sigma^+(u^n),e(\psi^u)\right) +\left(\sigma^-(u^n),e(\psi^u)\right) +((\Tilde{\varphi}^n)^2p^n, \operatorname{div} \psi^u)\\
    &\quad  + (1-\kappa)(\varphi^n\sigma^+(u^n) : e(u^n), \psi^\varphi - \varphi^n) +2(\varphi^n p^n \operatorname{div} u^n, \psi^\varphi-\varphi^n) \\
    &\quad + G_C\left(\frac{1}{\eps}(1-\varphi^n,\psi^\varphi-\varphi^n) + \eps(\nabla\varphi^n,\nabla (\psi^\varphi-\varphi^n))\right),
\end{align*}
for a given pressure $p^n\in L^\infty(\Omega)$. 
\end{Problem}

\subsection{Derivation of a complementarity system}
To treat the inequality in Problem \ref{problem:euler_lagrange_linearized}, 
we employ a Lagrange multiplier \RA{$\lambda^n$} by following \cite{kikuchi1988contact}. \RA{Initially, this Lagrange parameter exists in the dual space of $H^1(\Omega)$, but in \cite{Troeltzsch2004lagrange}, the author states that $\lambda^{\RA{n}} \in L^2(\Omega)$ for single constraint variational inequality problems with solution variables in $L^2(\Omega)$.} We assume that the argumentation can be transferred to Problem \ref{problem:euler_lagrange_linearized} and a rigorous proof is a goal for future work. With this, we can define
\begin{equation*}
    L^2_-(\Omega) \coloneqq \left\{v \in L^2(\Omega) \,|\, v\leq 0\quad \text{ a.e.\ in } \Omega \right\},
\end{equation*}
and
\begin{equation*}
    \mathcal{N}_+\coloneqq \left\{\mu \in L^2(\Omega) \,|\,  ( \mu, v )_{L^2(\Omega)}\leq 0\quad \forall v \in L^2_-(\Omega) \right\},
\end{equation*}
such that we can formulate a variational inequality system the classical $L^2$ 
inner products and $\mathcal{N}_+$ as solution set for the Lagrange multiplier $\lambda^n$.
\begin{Problem}\label{problem:euler_lagrange_with_lagrange_gelfand}
For a given $\varphi^0$ and for the incremental steps $t_n$ with $n=1,...,N$, find $U^n \in \mathcal{V} \times \mathcal{W}$ and $\lambda^n \in \mathcal{N}_+$ such that
\begin{align*}
    A(U^n)(\Psi) + (\lambda^n, \psi^\varphi) &= 0\quad \forall \Psi \in \mathcal{V} \times \mathcal{W} \cap L^\infty,\\
    (\lambda^n - \xi, \varphi^n - \varphi^{n-1}) &\geq 0 \quad \forall \xi \in \mathcal{N}_+,
\end{align*}
where $A(\cdot)(\cdot)$ is defined as before in Problem \ref{problem:euler_lagrange_linearized}.
\end{Problem}
The numerical method, which is introduced in Section \ref{sec:numerical_solution} to solve the inequality system, is designed to treat the variational inequality
\begin{equation*}
    (\lambda^n-\xi, \varphi^n - \varphi^{n-1}) \geq 0\quad \forall \xi \in \mathcal{N}_+,
\end{equation*}
in a \RA{complementarity formulation}. The following result states the equivalence of the 
previous variational inequality and a \RA{complementarity condition}.
\begin{Lemma}
\label{lemma_2_8}
The variational inequality
\begin{equation}\label{eq:var_ineq_lemma}
    (\lambda^n-\xi, \varphi^n - \varphi^{n-1}) \geq 0\quad \forall \xi \in \mathcal{N}_+,
\end{equation}
can equivalently be formulated \RA{as a complementarity condition} of the form
\begin{equation}\label{eq:complementarity_c}
    C(\varphi^n, \lambda^n)\coloneqq \lambda^n -\max\{0,\, \lambda^n +c(\varphi^n-\varphi^{n-1})\} =0,
\end{equation}
for every $c>0$ and the $\max$ operation defined as 
\begin{equation*}
    \max\{0,\, \lambda^n +c(\varphi^n-\varphi^{n-1})\} =\begin{cases}0 &\text{ in }\ \mathcal{I}, \\\lambda +c(\varphi^n-\varphi^{n-1}) &\text{ in }\ \mathcal{A},\end{cases}
\end{equation*}
where the inactive set $\mathcal{I} \subset \Omega$ and, the active set $\mathcal{A} \subset \Omega$
are defined such that  
\begin{align*}
    \lambda^n +c(\varphi^n-\varphi^{n-1}) &\leq 0 \text{ a.e.\ in } \mathcal{I},\\
    \lambda^n +c(\varphi^n-\varphi^{n-1}) &> 0 \text{ a.e.\ in } \mathcal{A}.
\end{align*}
These two situations can be explained as follows. For $\lambda^n = 0$ we are in the so-called inactive set, namely we solve the PDE part of phase-field. In $\mathcal{A}$ the constraint 
is active, `we sit on the obstacle', and we deal with $\lambda^n >0$.
Note, that $\mathcal{A}$ and $\mathcal{I}$ do not need to be connected, but note that $\Omega\setminus\mathcal{I}\cup\mathcal{A}$ is a null set. The sets can be understood as unions of all nonempty subsets of $\Omega$ with positive Lebesgue measure, on which the relations $>$ or $\leq$ are fulfilled almost everywhere.
\end{Lemma}
\begin{proof}
\RA{
The proof contains two major steps. Firstly, we prove the equivalence of \eqref{eq:var_ineq_lemma} and a strong formulation of the form 
\begin{align}
    \varphi^n - \varphi^{n-1} &\leq 0\quad \text{ a.e.\ in } \Omega,\label{eq:strong_ineq_1}\\
    \lambda^n &\geq 0\quad \text{ a.e.\ in } \Omega,\label{eq:strong_ineq_2}\\
    (\lambda^n, \varphi^n - \varphi^{n-1}) &= 0\label{eq:strong_ineq_3}.
\end{align}
Secondly, we show that \eqref{eq:strong_ineq_1}-\eqref{eq:strong_ineq_3} can be equivalently formulated as \eqref{eq:complementarity_c}.
}

\RA{\textit{Step 1:} ($\Rightarrow$)} Assume that \eqref{eq:var_ineq_lemma} is satisfied. For \eqref{eq:strong_ineq_1}, we assume the existence of a subset $\mathcal{T} \subset \Omega$ with a positive Lebesgue measure such that it holds $\varphi^n - \varphi^{n-1} > 0$ a.e.\ in $\mathcal{T}$. We define 
\begin{equation}
    \chi \coloneqq \begin{cases}2\lambda^n+ 1 &\text{ in } \mathcal{T}, \\
                                \lambda^n &\text{ in } \Omega \setminus \mathcal{T}. \end{cases}
\end{equation}
We have $\chi \in \mathcal{N}_+$ and obtain that the variational inequality must hold for $\mu = \chi$. But it holds
\begin{align*}
    (\lambda^n-\chi, \varphi^n-\varphi^{n-1}) &= (\lambda^n-\chi, \varphi^n-\varphi^{n-1})_{L^2(\mathcal{T})} + (\lambda^n-\chi,\varphi^n-\varphi^{n-1})_{L^2(\Omega \setminus \mathcal{T})} \\
                                              &= (\lambda^n-2\lambda^n \RA{-} 1, \varphi^n-\varphi^{n-1})_{L^2(\mathcal{T})} + (\lambda^n - \lambda^n,\varphi^n-\varphi^{n-1})_{L^2(\Omega \setminus \mathcal{T})}\\
                                              &= (-\lambda^n- 1, \varphi^n-\varphi^{n-1})_{L^2(\mathcal{T})}\\
                                              &= (\lambda^n+ 1, \varphi^{n-1}-\varphi^{n})_{L^2(\mathcal{T})} < 0,
\end{align*}
since $\varphi^{n-1}-\varphi^{n}<0$ a.e.\ in $\mathcal{T}$,  $\lambda^{\RA{n}} \in \mathcal{N}_+$ and $\lambda^n+1>0$ a.e.\ in $\Omega$. Summarizing, we obtain
\begin{equation*}
    (\lambda^n-\chi, \varphi^n-\varphi^{n-1}) < 0,
\end{equation*}
which is a contradiction to the assumption that \eqref{eq:var_ineq_lemma} is fulfilled. Thus, it must hold $\varphi^n-\varphi^{n-1} \leq 0$ a.e.\ in $\Omega$. For \eqref{eq:strong_ineq_2}, we assume $\lambda^n < 0$ a.e.\ in a subset $\mathcal{O} \subset \Omega$ with a positive Lebesgue measure. Then, we define 
\begin{equation*}
   v \coloneqq \begin{cases} \lambda^n    &\text{ in } \mathcal{O}\auth{,}\\
                                0       &\text{ in } \Omega \setminus \mathcal{O}\auth{.}\end{cases}
\end{equation*} 
Then, it holds $v \leq 0$ a.e.\ in $\Omega$. Thus, per definition of $\mathcal{N}_+$ and since $\lambda^n\in \mathcal{N}_+$, it must hold that $(\lambda^n, v) \leq 0$. But we find
\begin{equation*}
    (\lambda^{\RA{n}}, v) = \integral{\Omega}{}{\lambda^n v}{x} = \integral{\mathcal{O}}{}{\lambda^n v}{x} + \integral{\Omega\setminus\mathcal{O}}{}{\lambda^n v}{x} = \integral{\mathcal{O}}{}{\lambda^n \lambda^n}{x} + \integral{\Omega\setminus\mathcal{O}}{}{\lambda^n \cdot0}{x} = \integral{\mathcal{O}}{}{(\lambda^n)^2}{x} > 0,
\end{equation*}
which is a contradiction. Thus, it must hold $\lambda^n \geq 0$ a.e.~in $\Omega$. \\
Lastly, we derive \eqref{eq:strong_ineq_3}. We choose $\xi = 0$. Then, we have $\xi \in \mathcal{N}_+$ We obtain
\begin{equation*}
    (\lambda^n-0, \varphi^n - \varphi^{n-1}) = (\lambda^n, \varphi^n - \varphi^{n-1}) \geq 0.
\end{equation*}
In a similar way, we can set $\xi = 2\lambda^n$ to obtain 
\begin{equation*}
    (\lambda^n-2\lambda^n, \varphi^n - \varphi^{n-1}) = (-\lambda^n, \varphi^n - \varphi^{n-1})= -(\lambda^n, \varphi^n - \varphi^{n-1}) \geq 0,
\end{equation*}
which yields $(\lambda^n, \varphi^n-\varphi^{n-1}) \leq 0$. Combining both inequalities finally leads to $(\lambda^n, \varphi^n - \varphi^{n-1}) = 0$. \\
\RA{($\Leftarrow$)} Now, we assume that \eqref{eq:strong_ineq_1}-\eqref{eq:strong_ineq_3} hold true. Firstly, \eqref{eq:strong_ineq_2} validates the choice of $\mathcal{N}_+$ as solution space for $\lambda^n$. Let $v \in L^2_-$ be arbitrary. Then, we obtain $(\lambda^n, v) \leq 0$, and consequently $\lambda^n \in \mathcal{N}_+$. Now, let $\mu \ \in N_-$ be arbitrary. Due to \eqref{eq:strong_ineq_1}, we have $(\mu,\varphi^n-\varphi^{n-1}) \leq 0$, and thus $(-\mu,\varphi^n-\varphi^{n-1}) \geq 0$.
In combination with \eqref{eq:strong_ineq_3} we find 
\begin{equation*}
    (\lambda^n-\mu,\varphi^n-\varphi^{n-1}) \geq 0.
\end{equation*}
Since $\mu$ was chosen arbitrarily, we obtain \eqref{eq:var_ineq_lemma}.

\RA{\textit{Step 2:} ($\Rightarrow$)} Let \eqref{eq:strong_ineq_1}-\eqref{eq:strong_ineq_3} be fulfilled. We define $A,\, B \subset \Omega$ such that $\varphi^n -\varphi^{n-1} < 0$ a.e.\ in $A$ and $\lambda^n >0$ a.e.\ in $B$. As before, the sets $A$ and $B$ can be understood unions of all subsets on $\Omega$ with positive Lebesgue measure, on which the relations $>$ or $\leq$ are fulfilled almost everywhere. We start by proving that $A \cap B$ is a null set. For this, we assume that $A\cap B$ has a positive Lebesgue measure. We find $\lambda^n > 0$ a.e.\ in $A\cap B$, and $\varphi^n-\varphi^{n-1} <0$ a.e.\ in $A\cap B$. This yields
\begin{equation*}
    (\lambda^n, \varphi^n - \varphi^{n-1})_\Omega = (\lambda^n, \varphi^n - \varphi^{n-1})_{A\cap B} + (\lambda^n, \varphi^n - \varphi^{n-1})_{\Omega\setminus A\cap B} =  (\lambda^n, \varphi^n - \varphi^{n-1})_{A\cap B} <0,
\end{equation*}
which is a contradiction to \RA{\eqref{eq:strong_ineq_3}}. Thus, $A\cap B$ is a null set. In a next step, we prove that $C(\varphi^n, \lambda^n) =0$ a.e.\ in $\Omega$. It suffices to show that 
\begin{equation*}
    \restr{C(\varphi^n, \lambda^n)}{\mathcal{A}} = 0\quad \text{ a.e.\ in } \mathcal{A}, \quad \restr{C(\varphi^n, \lambda^n)}{\mathcal{I}} = 0\quad \text{ a.e.\ in } \mathcal{I},
\end{equation*}
hold true. \RA{I}n $\mathcal{A}$, we have 
\begin{equation*}
    \lambda^n \geq \lambda^n + c(\varphi^n-\varphi^{n-1}) > 0\quad \text{ a.e.\ in } \mathcal{A},
\end{equation*}
almost everywhere, which immediately yields $\lambda^n >0$ a.e.\ in $\mathcal{A}$. Due to previous findings, we consequently have $\varphi^n-\varphi^{n-1} = 0$ a.e.\ in $\mathcal{A}$, and obtain
\begin{equation*}
    \restr{C(\varphi^n, \lambda^n)}{\mathcal{A}}   = \lambda^n -\max\{0,\, \lambda^n +c(\varphi^n-\varphi^{n-1})\}= \lambda^n -\lambda^n -c(\varphi^n-\varphi^{n-1})= -c(\varphi^n - \varphi^{n-1})=0,
\end{equation*}
a.e.\ in $\mathcal{A}$. On $\mathcal{I}$ we have
\begin{equation*}
    c(\varphi^n-\varphi^{n-1})\leq \lambda^n + c(\varphi^n-\varphi^{n-1}) \leq 0\quad \text{ a.e.\ in } \mathcal{I},
\end{equation*}
which yields $\varphi^n-\varphi^{n-1} \leq 0$ a.e.\ in $\mathcal{I}$. With the same argumentation as before, we conclude $\lambda^n =0$ a.e.\ in $\mathcal{I}$ and obtain
\begin{equation*}
    \restr{C(\varphi^n, \lambda^n)}{\mathcal{I}}  = \lambda^n -\max\{0,\, \lambda^n +c(\varphi^n-\varphi^{n-1})\}= \lambda^n - 0=0,
\end{equation*}
a.e.\ in $\mathcal{I}$. Summarizing, \eqref{eq:complementarity_c} is fulfilled.\\
\RA{($\Leftarrow$)} Let \eqref{eq:complementarity_c} be fulfilled. As before, we have 
\begin{equation}\label{eq:lambda_larger_zero}
    \lambda^n + c(\varphi^n-\varphi^{n-1}) > 0\quad \text{ a.e.\ in }\mathcal{A},
\end{equation}
and
\begin{equation*}
    C(\varphi^n,\lambda^n)  = \lambda^n -\lambda^n - c(\varphi^n-\varphi^{n-1})=-c(\varphi^n - \varphi^{n-1})= 0,
\end{equation*}
a.e.\ in $\mathcal{A}$, which is equivalent to $\varphi^n - \varphi^{n-1} = 0$ a.e.\ in $\mathcal{A}$. Inserting this into \eqref{eq:lambda_larger_zero} yields $\lambda^n > 0 $ a.e.\ in $\mathcal{A}$. On $\mathcal{I}$ we have $\lambda^n +c(\varphi^n-\varphi^{n-1}) \leq 0 \text{ a.e.\ in } \mathcal{I}$,
and thus, \eqref{eq:complementarity_c} yields $\lambda^n = 0$ a.e.\ in $\mathcal{I}$. Thus, we find $\varphi^n -\varphi^{n-1} \leq 0$ a.e.\ in $\Omega$. To obtain \eqref{eq:strong_ineq_3}, we observe
\begin{equation*}
    (\lambda^n, \varphi^n-\varphi^{n-1}) = (\lambda^n, \varphi^n-\varphi^{n-1})_{\mathcal{A}} + (\lambda^n, \varphi^n-\varphi^{n-1})_{\mathcal{I}} = (\lambda^n, 0)_{\mathcal{A}} + (0, \varphi^n-\varphi^{n-1})_{\mathcal{I}} = 0.
\end{equation*}
With this, the proof is finished.
\end{proof}
In the following, the final system is stated and this is the starting point for 
the primal-dual active regularization.
\begin{Problem}\label{problem:euler_lagrange_with_complementarity}(Variational system with complementarity condition)
Given $\varphi^0$ and for the incremental steps $t_n$ with $n=1,...,N$, find $U^n =\{u^n,\varphi^n \}\in \mathcal{V} \times \mathcal{W}$ and $\lambda^n \in \mathcal{N}_+$ such that
\begin{align*}
    A(U^n)(\Psi) + (\lambda^n, \psi^\varphi) &= 0\quad \forall \Psi = \{\psi^u,\psi^{\varphi} \}\in \mathcal{V} \times \mathcal{W}\cap L^\infty,\\
    C(\varphi^n,\lambda^n) &= 0\quad \text{a.e.\ in }\Omega\auth{,}
\end{align*}
with $A(U^n)(\Psi)$ and $C(\varphi^n, \lambda^n)$ defined as before.
\end{Problem}

\section{A modified combined nonlinear Newton-type algorithm}
\label{sec:numerical_solution}
In this section, first, we recapitulate the principle numerical methods, which are used to solve the phase-field fracture problem introduced in Section \ref{sec:problem_formulation} and implemented in \texttt{pfm-cracks}~\cite{heister2020pfm-cracks} based on the finite element library \texttt{deal.II} \cite{dealii2021version94, dealii2019design}. 
The major computational features are 
an MPI parallelization (scalability tested on $1024$ cores \cite{HeiWi18_pamm}),
algebraic block-preconditiong using Trilinos \cite{trilinos1} inside a GMRES (generalized 
minimal residuals) linear iterative solver, a primal-dual 
active method for treating inequality constraints,
and predictor-corrector mesh adaptivity by choosing a given small regularization parameter 
\RA{$\varepsilon$} while guaranteeing $h<\varepsilon$, where $h$ is the local mesh size.

In the following two sections, we perform two improvements. The first modification enhances the performance of the primal-dual active set method 
by adjusting the constant $c$ introduced in \eqref{eq:complementarity_c}.
With the second modification, we improve the accuracy in time by updating the linearization term $\Tilde{\varphi}$ until convergence (into the monolithic limit) within one incremental step. The resulting final algorithm is then 
implemented in \texttt{pfm-cracks} and constitutes an extension to 
existing published work.

\subsection{Newton's method}
In each incremental step, the first equation of Problem \ref{problem:euler_lagrange_with_complementarity} can be solved via a Newton method. Considering incremental step $n$, we seek for $U^n \coloneqq \{u^n,\varphi^n\}\in \mathcal{V} \times \mathcal{W}$ and $\lambda^n \in \mathcal{N}_+$. The solution $U^n$ is found via iterating over $k=1,2,3,\ldots$ until convergence, such that
\begin{equation*}
    A'(U^{n,k})(\delta U^{n,k+1}, \Psi) + (\lambda^{n,k},\psi^\varphi) = -A(U^{n,k})(\Psi) \quad \forall \Psi \in \mathcal{V}\times \mathcal{W}
\end{equation*}
with respect to
\begin{equation*}
    C(\varphi^{n,k} + \delta \varphi^{n,k+1},\lambda^{n,k}) =0\quad \text{ a.e.\ in } \Omega,
\end{equation*}
for the update $\delta U^{n,k+1}$ and the Lagrange multiplier $\lambda^{n,k+1}$ and updating via 
\begin{equation*}
    U^{n,k+1} = U^{n,k} + \delta U^{n,k+1}.
\end{equation*}
The Jacobian $ A'(U^{n,k})(\delta U^{n,k+1}, \Phi)$ is given by
\begin{align*}
    A'(U^{n,k})(\delta U^{n,k+1}, \Psi) &= \left(g(\Tilde{\varphi}^{n}) \sigma^+(\delta u^{n,k+1}), e(\psi^u) \right) + \left(\sigma^-(\delta u^{n,k+1}), e(\psi^u)\right) \\
                                &\qquad + (1-\kappa)\left(\delta \varphi^{n,k+1} \sigma^+(u^{n,k}) : e(u^{n,k}) + 2 \varphi^{n,k} \sigma^+(\delta u^{n,k+1}) : e(u^{n,k}), \psi^\varphi\right) \\
                                &\qquad + 2p\left(\delta\varphi^{n,k+1} \operatorname{div} u^{n,k} + \varphi^{n,k} \operatorname{div} \delta u^{n,k+1}, \psi^\varphi\right) \\
                                &\qquad + G_C \left( \frac{1}{\eps}(\delta \varphi^{n,k\RA{+1}}, \psi^\varphi) + \eps (\nabla \delta \varphi^{n,k+1}, \nabla \psi^\varphi)\right),
\end{align*}
and $A(U^{n,k})(\Phi)$ is defined as before. To treat the complementarity condition, we introduce the primal-dual active set method.  

\subsection{The primal-dual active set method}\label{subsec:classic_active_set}
The primal-dual active set method \RA{(PDAS)}, introduced for constrained inequality systems in \cite{bergounioux2000comparison, bergounioux1999primal, ito2000augmented, ito2000optimal},
and shown under certain assumptions to be a semi-smooth Newton method \cite{hintermueller2002semismooth}, and applied on the phase-field fracture model in \cite{heister2015primal}, is based on considerations made in Section \ref{sec:problem_formulation}. The idea is to split the domain into two subdomains \auth{in each incremental step $n$}. On one subdomain, the inactive set $\mathcal{I}^{\auth{n}}$, the inequality constraint $\varphi^n -\varphi^{n-1} \leq 0$ is fulfilled strictly. On the other subdomain, the active set $\mathcal{A}^{\auth{n}}$, it holds $\varphi^n -\varphi^{n-1} = 0$. A priori, these two sets are not known but combined with the previously introduced Newton method, we obtain a prediction algorithm. Based on the complementarity condition, \auth{i.e., \eqref{eq:complementarity_c}}
\begin{equation*}
    C(\varphi^{n,k}+\delta \varphi^{n,k+1},\lambda^{n,k}) = 0,
\end{equation*}
we want to determine the active set $\mathcal{A}^{\auth{n,k}}$ and the inactive set $\mathcal{I}^{\auth{n,k}}$ such that 
\begin{equation}
\label{eq_determine_active_set}
    \lambda^{n\auth{,k}} + c(\varphi^{n\auth{,k}}-\varphi^{n-1}) > 0\quad \text{a.e.\ in } \mathcal{A}^{\auth{n,k}}, \quad  \lambda^{n\auth{,k}} + c(\varphi^{n\auth{,k}}-\varphi^{n-1}) \leq 0\quad \text{a.e.\ in } \mathcal{I}^{\auth{n,k}},
\end{equation}
\auth{and iterate until the active set does not change within two consecutive Newton iterations.} As before, $\mathcal{A}^{\auth{n}}$ can be understood as the union of all subsets of $\Omega$ on which $\lambda^n + c(\varphi^n-\varphi^{n-1}) > 0$ is fulfilled almost everywhere. The inactive set $\mathcal{I}^{\auth{n}}$ is then the complement of $\mathcal{A}^{\auth{n}}$ with respect to $\Omega$. If $\mathcal{A}^{\auth{n}}$ and $\mathcal{I}^{\auth{n}}$ are known, we can set $\lambda^n=0$ on $\mathcal{I}^{\auth{n}}$, see Lemma \ref{lemma_2_8}, 
and treat the problem as an unconstrained problem. On $\mathcal{A}^{\auth{n}}$, we set $\varphi^n = \varphi^{n-1}$ and there is nothing to do (for the phase-field). The resulting 
scheme is given in Algorithm \ref{alg:active_set_cont}.
\begin{algorithm}
\caption{(Primal-dual active set method)}\label{alg:active_set_cont}
\begin{algorithmic}[1]
\State \RA{Set iteration index $k=0$}
\RA{\While{$\mathcal{A}^{n,k} \neq \mathcal{A}^{n,k+1}$}
\State Determine the active set $\mathcal{A}^{n,k}$ and inactive set $\mathcal{I}^{n,k}$ with \eqref{eq_determine_active_set}
\State Find $\delta U^{n,k+1} \in \mathcal{V}\times \mathcal{W}$ and $\lambda^{n,k+1}\in \mathcal{N}_+$ with solving
\begin{align*}
    A'(U^{n,k})(\delta U^{n,k+1}, \Phi) + (\lambda^{n,k+1},\psi) &= -A(U^{n,k})(\Phi),\quad \forall \Phi\coloneqq \{v,\psi\} \in \mathcal{V}\times \mathcal{W},\\
    \delta \varphi^{n,k+1} &= 0\quad \text{ on } \mathcal{A}^k,\\
    \lambda^{n,k+1} &= 0\quad \text{ on } \mathcal{I}^k.
 \end{align*}
\State Update the solution to obtain $U^{n,k+1}$ via
\begin{equation*}
    U^{n,k+1} = U^{n,k} + \delta U^{n,k+1}.
\end{equation*}
\State \RA{Update iteration index $k=k+1$}
\EndWhile}
\end{algorithmic}
\end{algorithm}

\subsection{Discretization}
We discretize Problem \ref{problem:euler_lagrange_with_complementarity} using a finite element method with bilinear (2d) or trilinear (3d) elements $Q_c^1$ \cite{ciarlet2002fem} 
for both, the displacement function and the phase-field function. The discrete function spaces read
\begin{align*}
    \mathcal{V}_h &\coloneqq \bigg\{ u_h \in \mathcal{V}, \, \restr{u_h}{K} \in \left[Q_1^c(K)\right]^d,\quad \forall\,  K \in \mathcal{T}_h \bigg \}, \\
    \mathcal{W}_h &\coloneqq \bigg\{ \varphi_h \in \mathcal{W}, \, \restr{\varphi_h}{K} \in Q_1^c(K),\quad \forall \, K \in \mathcal{T}_h \bigg\},\\
    \mathcal{N}_h &\coloneqq \bigg\{ \lambda_h \in \mathcal{N}_+, \, \restr{\lambda_h}{K} \in Q_1^c(K),\quad \forall \, K \in \mathcal{T}_h \bigg\},
\end{align*}
where $K\in\mathcal{T}_h$ is the finite element \auth{and} 
$\mathcal{T}_h$ denotes the 
decomposition of the domain $\Omega$ into a mesh, e.g., \cite{ciarlet2002fem}. 
The spatially discretized system is 
\auth{formulated as follows:}
\begin{Problem}\label{problem:discretized_euler_lagrange_with_complementarity}(Discretized system with complementarity condition)
Given $\varphi_h^0$ and for the incremental steps $t_n$ with $n=1,...,N$, find $(u_h^n,\varphi_h^n, \lambda_h^n) \in \mathcal{V}_h \times \mathcal{W}_h \times \mathcal{N}_h$ such that
\begin{align*}
    A(u_h^n,\varphi_h^n)(\Phi_h) + (\lambda_h^n, \psi_h) &= 0\quad \forall \Phi_h\coloneqq(v_h,\psi_h) \in \mathcal{V}_h \times \mathcal{W}_h,\\
    C(\varphi_h^n,\lambda_h^n) &= 0\quad \forall x \in \Omega,
\end{align*}
with $A(u_h^n,\varphi_h^n)(v_h,\psi_h)$ and $C(\varphi_h^n, \lambda_h^n)$ 
defined as before in Lemma \ref{lemma_2_8} with the point-wise maximum operation. 
\end{Problem}
The application of a Galerkin ansatz with primitive ansatz and test functions of the form
\begin{align*}
    \Phi_{h,i} &= \begin{bmatrix}\chi_i^u\\0\end{bmatrix} \text{ for } i=1,...,N_u,\\
    \Phi_{h,N_u+i} &= \begin{bmatrix}0\\\chi_i^\varphi\end{bmatrix} \text{ for } i=1,...,N_\varphi,
\end{align*}
 where $N_u$ is the number of degrees of freedom in $u_h$ and $N_\varphi$ is the number of degrees of freedom in $\varphi_h$, leads to a system of the form
\begin{equation*}
    \begin{bmatrix}
    M & B\\
    B^T & 0
    \end{bmatrix}
    \begin{bmatrix}
    \delta U_h^{n,k+1}\\
    \lambda_h^{k+1} 
    \end{bmatrix}
    =\begin{bmatrix}
    F\\
    0
    \end{bmatrix},
\end{equation*}
where $B$ is a mass matrix and $M$ and $F$ are given by 
\begin{equation*}
M = \left[\begin{array}{cc}
     M^{uu} & M^{u\varphi}\\
     M^{\varphi u} & M^{\varphi \varphi}
\end{array} \right],
\quad F = \left[\begin{array}{c}
     F^{u}   \\
     F^{\varphi}
\end{array} \right],
\end{equation*}
with the block entries
\begin{align*}
    M_{ij}^{uu} &= \left(g(\Tilde{\varphi}^n) \sigma^+(\chi_j^u), e(\chi_i^u) \right) + \left(\sigma^-(\chi_j^u), e(\chi_i^u)\right),\\
    M_{ij}^{\varphi u} &= 2(1-\kappa)\left(\varphi_h^{n,k} \sigma^+(\chi_j^u) : e(u_h^{n,k}), \chi_i^\varphi\right)+ 2p\left(\varphi_h^{n,k} \operatorname{div} (\chi_j^u), \chi_i^\varphi\right), \\
    M_{ij}^{u\varphi} &= 0, \\
    M_{ij}^{\varphi \varphi} &= (1-\kappa)\left(\sigma^+(u_h^{n,k}) : e(u_h^{n,k})\chi_j^\varphi, \chi_i^\varphi\right)+ 2p\left(\operatorname{div} (u_h^{n,k}) \chi_j^\varphi, \chi_i^\varphi\right) \\
        &\qquad + G_C \left(\frac{1}{\eps}\left(\chi_j^\varphi, \chi_i^\varphi\right) + \eps \left(\nabla \chi_j^\varphi, \nabla \chi_i^\varphi\right)\right),\\
&\text{and}\\
\begin{split}
    F_{ij}^{u} &= -A(U_h^{n,k})(\chi_i^u) = -\left(\left[(1-\kappa)(\Tilde{\varphi}_h^n)^2+\kappa\right]\sigma^+(u_h^{n,k}),e(\chi_i^u)\right) - \left(\sigma^-(u_h^{n,k}),e(\chi_i^u)\right)\\
    &\qquad  \qquad \qquad \qquad \quad - \left((\Tilde{\varphi}_h^n)^2p, \operatorname{div} (\chi_i^u)\right),
\end{split}
\\[1ex]
\begin{split}
    F_{ij}^{\varphi} &= -A(U_h^{n,k})(\chi_i^\varphi) = -(1-\kappa)\left(\varphi_h^{n,k} \sigma^+(u_h^{n,k}) : e(u_h^{n,k}), \chi_i^\varphi\right) - 2 \left(\varphi_h^{n,k} p \operatorname{div} (u_h^{n,k}), \chi_i^\varphi\right)\\
    &\qquad  \qquad \qquad \qquad \quad- G_C\left(\frac{1}{\eps}\left(1-\varphi_h^{n,k},\chi_i^\varphi\right) + \eps\left(\nabla\varphi_h^{n,k},\nabla \chi_i^\varphi\right)\right).
    \end{split}
\end{align*}
\RA{We notice that the block 
$M_{ij}^{u\varphi}$ is zero due to the previously applied linearization in $\varphi$ in the 
displacement equation in Problem \ref{problem:euler_lagrange_linearized}. Consequently, in the Newton system matrix, the corresponding directional derivative vanishes and it holds $M_{ij}^{u\varphi}=0$. The main purpose 
is a robust nonlinear and linear solution; see also Remark \ref{remark:linear_solver}.
}

Since the $u_h^{n,k}, \, \varphi_h^{n,k},\, \lambda_h^{n,k}$ are element-wise of polynomial structure, we can compute the active and inactive set point-wise:
\begin{align*}
    \mathcal{A}^{\auth{n,}k} &= \left\{ x \, | \, \lambda_h^{n,k}(x)+c(\varphi_h^{n,k}(x) -\varphi_h^{n-1}(x))>0 \right\},\\
    \mathcal{I}^{\auth{n,}k} &= \left\{ x \, | \, \lambda_h^{n,k}(x)+c(\varphi_h^{n,k}(x) -\varphi_h^{n-1}(x))\leq0 \right\}.
\end{align*}
Given the fully discretized system, we can formulate the primal-dual active set method as it is implemented in \texttt{pfm-cracks} in Algorithm \ref{alg:active_set_disc}.

\begin{algorithm}[h!]
\caption{(Primal-dual active set method \RA{with backtracking line search})}\label{alg:active_set_disc}
\begin{algorithmic}[1]
\State \RA{Set iteration index $k=0$}
\RA{\While{$\left(\mathcal{A}^{n,k-1} \neq \mathcal{A}^{n,k}\right)$ \textbf{or} $\left(\Tilde{R}(U_h^{n,k}) > \text{TOL}_N\right)$}
    \State Assemble the residual $R(U_h^{n,k})$
    \State Compute the active set $\mathcal{A}^{n,k} = \left\{ x_i  \bigg| \left[B\right]_{ii}^{-1} \left[R(U_h^{n,k})\right]_i + c(\varphi_{h,i}^{n,k} - \varphi_{h,i}^{n-1}) > 0\right\} $
    \State Set $\varphi_h^{n,k} = \varphi_h^{n-1}$ on $\mathcal{A}^{n,k}$
    \State Assemble the system matrix $M$ (Newton Jacobian) and the right-hand side $F=R(U_h^{n,k})$
    \State Eliminate rows/columns in $\mathcal{A}^{n,k}$ from $M$ and $F$ to obtain $\Tilde{M}$ and $\Tilde{F}=\tilde{R}(U_h^{n,k})$ 
    \State Solve the linear system $\Tilde{M}\delta U_h^{n,k+1} = \Tilde{F}$ with GMRES and AMG preconditioner \cite{heister2020pfm-cracks}
    \State Choose maximum number of line search iterations $l_{\max}$
    \State Choose line search damping parameter $0 < \omega \leq 1 $
    \For{$l=1:l_{\max}$}
        \State Update the solution with $U_h^{n,k+1} = U_h^{n,k} + \delta U_h^{n,k+1}$
        \State Assemble the new residual $\tilde{R}(U_h^{n,k+1})$
        \If{$\lVert \tilde{R}(U_h^{n,k+1}) \rVert_2 < \lVert \tilde{R}(U_h^{n,k}) \rVert_2$}
            \State \textbf{break}
        \Else
            \State Adjust the Newton update with $\delta U_h^{n,k+1} := \omega^l\delta U_h^{n,k+1}$
        \EndIf
    \EndFor
    \State Update iteration index $k=k+1$
\EndWhile
}
\end{algorithmic}
\end{algorithm}

\RA{
\begin{Remark}
\label{remark:linear_solver}
In line 8, the reduced linear system is solved with a GMRES (generalized minimal residual) method \cite{saad1986gmres} and 
algebraic multigrid preconditioning (AMG) \cite{trilinos1}. In this work, the implementation as it is from \texttt{pfm-cracks} \cite{heister2020pfm-cracks} is utilized. In our numerical tests (Section \ref{sec:numerical_tests}), we observe 
in all simulations between $10-40$ linear iterations. This is in agreement with the results obtained in \cite{HeiWi18_pamm}[Table 1]. The main reason for the excellent performance is twofold. First, it is the 
triangular block structure of $\tilde M$ due to the zero block $M_{ij}^{u\varphi}$ as previously discussed. Second, the diagonal 
terms in $\tilde M$ are of elliptic type, which is well-known that multigrid methods perform very well.
\end{Remark}
}

\RA{\begin{Remark}
The lines 9-19 describe a classical backtracking line search algorithm, where the Newton update is damped with a damping parameter $\omega\in (0,1]$, if the updated solution $U_h^{n,k+1} = U_h^{n,k} + \delta U_h^{n,k+1}$ does not reduce the residual norm. In all experiments in Section \ref{sec:numerical_tests}, we use $l_{\max} = 10$ and $\omega = 0.6$.
\end{Remark}}

\RA{\begin{Remark}
Note that we deal with two systems of equations in Algorithm \ref{alg:active_set_disc}: the global nonlinear system and the reduced linear system. The nonlinear system consists of the matrix $M$ and the right-hand-side $F$. The reduced linear system, defined by the matrix $\Tilde{M}$ and $\Tilde{F}$ only contains the equations of the nonlinear system, which belong to the inactive set. The residual of the full nonlinear system is then given by $R$ whereas $\Tilde{R}$ denotes the residual of the reduced linear system. From an implementation point of view, this is realized by setting constraints to the system such that we enforce the phase-field to remain the same on the active degrees of freedom.
\end{Remark}}

\subsection{Modified combined Newton active set algorithms}\label{subsec:mod_newton_active_set}
This section is dedicated to presenting adjustments to the primal-dual active set method based on an analysis of the active set constant $c>0$ (see again Lemma \ref{lemma_2_8} and 
Algorithm \ref{alg:active_set_cont}). In \cite{kaerkkaeinen2003augmented}, the authors prove for an obstacle problem that the primal-dual active set converges for any sufficiently large $c>0$. Furthermore, they point out that $c$ only influences the first active set iteration in theory. Similar observations were made in \cite{hueber2005primaldual,popp2009finite,schroeder2016semismooth}. The authors state on the one hand that a constant $c$ of magnitude around the Youngs modulus is reasonable from an engineer's perspective, while they also point out on the other hand, that different settings for $c$ do not affect the solution but only the algorithmic performance. These results are a motivation to further investigate the influence of $c$ for our primal-dual active set phase-field fracture formulation. 
We start by pointing out the bottleneck of the above described primal-dual active set algorithm. 
Specifically on fine meshes, we often run into convergence issues of the active set, whereas the residual converges comparably fast (see e.g., \cite{heister2015primal}[Fig. 14]). 

\subsubsection{Investigation of the active set constant \texorpdfstring{$c$}{c}}
In the following, our objective is to illustrate the influence of the active set constant $c$ and how to adjust it to reduce the number of active set iterations.
To this end, we begin with an investigation of the role of the active set constant $c$ in the algorithm. It is involved in the classification of the active set $\mathcal{A}^{k}$ as a degree of freedom $x_i$ is classified as active (in iteration $k$ \auth{within the current incremental step}), if
\begin{equation*}
    \lambda_{h,i}^k + c(\varphi_{h,i}^k - \varphi_{h,i}^{\operatorname{old}}) > 0,
\end{equation*}
\auth{where $\varphi_{h,i}^{\operatorname{old}}$ is the value of the phase field solution of the previous incremental step at degree of freedom $i$.}
We identify nine different situations, depending on the sign of $\lambda_{h,i}^k$ and $(\varphi_{h,i}^k-\varphi_{h,i}^{\operatorname{old}})$:
\begin{enumerate}
        \item $\varphi_{h,i}^k - \varphi_{h,i}^{\operatorname{old}} = 0, \text{ and } \lambda^k_{h,i} > 0 \Rightarrow \varphi_{h,i}^k= \varphi_{h,i}^{\operatorname{old}},\, \delta \varphi^{k+1}_{h,i} = 0 \text{ and }  \lambda^{k+1}_{h,i} = B_{ii}^{-1} F_i$
        \item $\varphi_{h,i}^k - \varphi_{h,i}^{\operatorname{old}} = 0, \text{ and } \lambda^k_{h,i} = 0 \Rightarrow \delta \varphi^{k+1}_{h,i} \text{ as solution of the system and }  \lambda^{k+1}_{h,i} = 0$
        \item $\varphi_{h,i}^k - \varphi_{h,i}^{\operatorname{old}} = 0, \text{ and } \lambda^k_{h,i} < 0 \Rightarrow \varphi_{h,i}^k= \varphi_{h,i}^{\operatorname{old}},\, \delta \varphi^{k+1}_{h,i} = 0 \text{ and }  \lambda^{k+1}_{h,i} = B_{ii}^{-1} F_i$
        \item $\varphi_{h,i}^k - \varphi_{h,i}^{\operatorname{old}} > 0, \text{ and } \lambda^k_{h,i} > 0 \Rightarrow \delta \varphi^{k+1}_{h,i} = 0 \text{ and }  \lambda^{k+1}_{h,i} = B_{ii}^{-1} F_i$
        \item $\varphi_{h,i}^k - \varphi_{h,i}^{\operatorname{old}} > 0, \text{ and } \lambda^k_{h,i} = 0 \Rightarrow \varphi_{h,i}^k= \varphi_{h,i}^{\operatorname{old}},\, \delta \varphi^{k+1}_{h,i} = 0 \text{ and }  \lambda^{k+1}_{h,i} = B_{ii}^{-1} F_i$
        \item $\varphi_{h,i}^k - \varphi_{h,i}^{\operatorname{old}} > 0, \text{ and } \lambda^k_{h,i} < 0 \Rightarrow \begin{cases} \varphi_{h,i}^k= \varphi_{h,i}^{\operatorname{old}},\, \delta \varphi_i^{k+1} = 0, \, \lambda^{k+1}_{h,i} = B_{ii}^{-1} F_i &\text{if } |\lambda_{h,i}^{k}| < c(\varphi_{h,i}^k-\varphi_{h,i}^{\operatorname{old}}),\\
                    \delta \varphi_{h,i}^{k+1} \text{ as solution}, \, \lambda^{k+1}_{h,i} = 0 &\text{otherwise.}\end{cases}$
        \item  $\varphi_{h,i}^k - \varphi_{h,i}^{\operatorname{old}} < 0, \text{ and } \lambda^k_{h,i} > 0\Rightarrow \begin{cases}  \delta \varphi_{h,i}^{k+1} \text{ as solution}, \, \lambda^{k+1}_{h,i} = 0 &\text{if }  \lambda_{h,i}^k \leq |c(\varphi_{h,i}^k-\varphi_{h,i}^{\operatorname{old}})|  \\
        \varphi_{h,i}^k= \varphi_{h,i}^{\operatorname{old}},\,\delta \varphi_{h,i}^{k+1} = 0, \, \lambda^{k+1}_{h,i} = B_{ii}^{-1} F_i    &\text{otherwise}\end{cases}$      
        \item $\varphi_{h,i}^k - \varphi_{h,i}^{\operatorname{old}} < 0, \text{ and } \lambda^k_{h,i} = 0 \Rightarrow  \delta \varphi_{h,i}^{k+1} \text{ as solution}, \, \lambda^{k+1}_{h,i} = 0$
        \item $\varphi_{h,i}^k - \varphi_{h,i}^{\operatorname{old}} < 0, \text{ and } \lambda^k_{h,i} < 0 \Rightarrow  \delta \varphi_{h,i}^{k+1} \text{ as solution}, \, \lambda^{k+1}_{h,i} = 0$.
\end{enumerate}
We observe that the active set constant $c$ only has an influence on the classification when $\lambda_{h,i}^k$ and $\varphi_{h,i}^k - \varphi_{h,i}^{\operatorname{old}}$ have different signs. In situation No. 6, $\varphi_{h,i}^k$ shows crack healing behaviour, thus \auth{we} do not want to accept it as a solution and set it to $\varphi_{h,i}^{\operatorname{old}}$. We achieve this, if 
\begin{equation*}
    |\lambda_{h,i}^{k}| < c(\varphi_{h,i}^k-\varphi_{h,i}^{\operatorname{old}}),
\end{equation*}
i.e.
\begin{equation*}
    \frac{|\lambda_{h,i}^{k}|}{(\varphi_{h,i}^k-\varphi_{h,i}^{\operatorname{old}})} < c.
\end{equation*}
In situation No. 7, we do not violate the constraint with the solution of the $k$th iteration, thus, we still want to classify this degree of freedom as inactive. This can be achieved via 
\begin{equation*}
     \lambda_{h,i}^k \leq c|(\varphi_{h,i}^k-\varphi_{h,i}^{\operatorname{old}})|,
\end{equation*}
i.e.
\begin{equation*}
     \frac{\lambda_{h,i}^k}{|(\varphi_{h,i}^k-\varphi_{h,i}^{\operatorname{old}})|} < c.
\end{equation*}
Summarizing, we can formulate a condition for $c$:
\begin{equation*}
     \left|\frac{\lambda_{h,i}^k}{(\varphi_{h,i}^k-\varphi_{h,i}^{\operatorname{old}})} \right|< c. 
\end{equation*}
Thus, any $c$ larger than the lower bound is sufficiently large. 
\subsubsection{Proposed adjustments and definition of four cases}\label{subsubsec:4_cases}
\RA{In the following, we propose four different cases for implementing the Newton active set algorithm. The basis for these cases is Algorithm \ref{alg:active_set_disc}. For better readability and since the concepts are the same in each incremental step, we drop the incremental index $n$. Both stopping criteria, i.e., if not stated otherwise, require the active set to converge and the residual-norm to fall below a certain tolerance.}
For our adjustment, we choose 
\begin{equation*}
     c = c^k \coloneqq 2\left|\frac{\lambda_{h,i}^k}{(\varphi_{h,i}^k-\varphi_{h,i}^{\operatorname{old}})} \right|, 
\end{equation*}
thus, in contrast to before, $c$ changes in every \auth{Newton} iteration. This could be avoided by iterating until convergence, saving the largest $c^k$ and then restarting the iteration. But this is an unnecessary computational cost and in our opinion, a varying $c^k$ does not lead to any conflicts. In Section \ref{sec:numerical_tests}, we will perform several experiments to observe the performance boost of this adjustment. Based on the previous findings, we suggest four different cases:
\begin{itemize}
\item \textbf{Case 1:} We iterate as long as the active set does not change within 2 iterations with a constant $c = 10E$, where $E$ is Young's modulus. Let $k$ be the iteration index and $\mathcal{A}^k$ the active set of iteration $k$, we stop, when $\mathcal{A}^k = \mathcal{A}^{k+1}$ \RA{and $\lVert \Tilde{R}(U_h^{k+1}) \rVert_2< \operatorname{TOL}_N$}.
\item \textbf{Case 2:} The classification of active/inactive set proceeds as in \textbf{Case 1}, but with the modified $c$ set as 
\begin{equation*}
     c = c^k = 2\left|\frac{\lambda_{h,i}^k}{(\varphi_{h,i}^k-\varphi_{h,i}^{\operatorname{old}})} \right|. 
\end{equation*}
\RA{Apart from this, everything is similar to Algorithm \ref{alg:active_set_disc} including the stopping criteria: we stop, when $\mathcal{A}^k = \mathcal{A}^{k+1}$ and $\lVert \Tilde{R}(U_h^{k+1}) \rVert_2< \operatorname{TOL}_N$}.
\item \textbf{Case 3:} \RA{The classification of the active/inactive set proceeds as in \textbf{Case 1} with $c=10E$. But in this case, we weaken the active set stopping criterion, i.e.\ we do not enforce $\mathcal{A}^k = \mathcal{A}^{k+1}$ for termination anymore. Instead, we perform as much Newton active set iterations as needed to achieve $\lVert \Tilde{R}(U_h^{k+1})\rVert_2 <~\operatorname{TOL}_N$. When this is fulfilled, we only perform $10$ more Newton active set iterations until we stop. 
The number $10$ is chosen heuristically based on our experiences 
in this paper. It is a compromise between sufficiently many iterations to guess that we may have 
converged and computational cost by not adding too many additional iterations.}
\item \textbf{Case 4:} \RA{The classification of the active/inactive set proceeds as in \textbf{Case 1} with $c=10E$. But in this case, we completely omit the active set stopping criteria. This means, we do not require $\mathcal{A}^k = \mathcal{A}^{k+1}$ but stop immediately as soon as $\lVert \Tilde{R}(U_h^{k+1})\rVert_2 < \operatorname{TOL}_N$ is reached.}
\end{itemize}

\section{Iteration on the linearization into the monolithic limit}
\label{sec:iteration_on_extra}
Both introduced linearization strategies (Section \ref{sec_PFF_linearization}) to treat the nonlinearity in the displacement equation come with limitations of the above model. 
They lead to so-called time-lagging behaviour (temporal/incremental discretization error), 
where the crack grows slower than physics of the governing model suggest; \auth{see for instance \cite{Wi17_SISC}[Fig. 3].}
Based on concepts developed in \cite{wick2020multiphysics}[Section 7.7.3], we employ an iteration on the linearization, such that the iteration converges into the monolithic limit. 
The idea is to employ an additional fixed-point iteration.
In this fashion, we iterate until the $L^2$-difference between two consecutive iterative solutions is smaller than a predefined tolerance. This concept can be applied to both the extrapolation (ItE) and the linearization by using the solution from the previous incremental step (ItOTS); see again 
Section \ref{sec_PFF_linearization}. 
For the former, numerical studies were already done in \cite{LAMPRON2021114091}. 

The main objective is to further investigate ItE in practice (since from a pure 
mathematical point of view the extrapolation is only heuristic since 
no regularity in time can be ensured since $\partial_t \varphi$ 
is only a bounded measure \cite{MiWheWi15b}[p. 1384, Theorem 1] shown for 
a decoupled formulation; but for extrapolation more regularity 
is required). For comparison, we utilize as well ItOTS.
Along with these iterations into the monolithic limit, 
a challenge is its combination with predictor-corrector adaptivity.
In the predictor-corrector algorithm \cite{heister2015primal}[Section 4]
first a new crack path is predicted, and then solved again on the new mesh.
Iterating on the new mesh into the monolithic limit may result into 
a previous under-estimation of the crack path, for which the predictor-corrector 
algorithm needs to be re-started, because the condition $h<\varepsilon$ might 
be violated since the predicted refinement area was not large enough. 
Of course, this double iterations can become 
quite expensive and finally it is a compromise between two typical 
numerical demands, namely efficiency and 
\auth{accuracy}.

Therefore, 
our recommendation is not that we must use the full algorithm in all
various cases, it is rather a decision choice what is more important: 
very small internal length scales $\varepsilon$ (thus predictor-corrector adaptivity
is likely needed), overall accuracy of the solution, overall efficiency
(including adaptivity and/or parallel computing). 

In consequence, our newly proposed final Algorithm \ref{alg:full_algorithm} includes the following iterations: The loop over the incremental steps, the iteration on the linearization, possible predictor-corrector adaptive mesh refinement iteration, the nonlinear Newton active set iteration and the linear solver iteration. 
\begin{algorithm}
\caption{(Quasi-monolithic solution algorithm including iteration on the linearization)}\label{alg:full_algorithm}
\begin{algorithmic}[1]
\State Setup the system                                 \Comment{initialize grid $\mathcal{T}_h$, parameters, etc.}
\State Set $\mathcal{T}_h^{\operatorname{old}} = \mathcal{T}_h$ 
\For{$n=1,2,...$}\Comment{timestep loop}
    \State Set $\operatorname{changed_{mesh}}=$\textbf{true} 
    \While{$(\operatorname{changed_{mesh}})$}                        \Comment{adaptive predictor-corrector refinement}
        \State Set $\varphi_h^{n,-1} = \varphi_h^{n-2}$ \Comment{preparations for iteration on the linearization}
        \State Set $\varphi_h^{n,0} = \varphi_h^{n-1}$  \Comment{preparations for iteration on the linearization}
        \State Set $j = 0$                              \Comment{index for iteration on the linearization}
        \While{$\left(\lVert \varphi_h^{n,j}-\varphi_h^{n,j-1}\rVert_2 \geq \text{TOL}_{\operatorname{ItL}}\right)$ \textbf{or} $ \left(j < 1\right)$} \Comment{iteration on the linearization}
            \State Compute the linearization $\tilde{\varphi}_h^{n,j}$ \Comment{compute the linearization (Section \ref{sec_PFF_linearization})}
            \State Set $j = j+1$ \Comment{update \auth{ItL} index}
            \State Solve system with PDAS to obtain $\{u_h^{n,j},\varphi_h^{n,j}\}$  \Comment{ Algorithm \ref{alg:active_set_disc} with linearization in assembly}
        \EndWhile
        \State Refine the mesh to obtain new mesh $\mathcal{T}_h^{\operatorname{new}}$ \Comment{predictor-corrector scheme \cite{heister2015primal}}
        \If{$\mathcal{T}_h^{\operatorname{new}} = \mathcal{T}_h^{\operatorname{old}}$} \Comment{Check whether mesh changed or not}
            \State Set $\operatorname{changed_{mesh}}=$ \textbf{false} \Comment{If did not change, leave while loop}
        \Else
            \State Set $\mathcal{T}_h^{\operatorname{old}} = \mathcal{T}_h^{\operatorname{new}}$ \Comment{Else, save current mesh and goto line 6}
        \EndIf
    \EndWhile
\EndFor
\end{algorithmic}
\end{algorithm}

\newpage
\section{Numerical examples}\label{sec:numerical_tests}
In \auth{this section}, we propose several numerical \auth{tests} for the previously introduced \auth{algorithms}. In the first two examples, the Sneddon test in two and three dimensions, we compare the required Newton active set iterations for solving the problem for the four different cases \auth{(Section \ref{subsubsec:4_cases})}. Since these tests are steady-state, i.e. we have a non-growing fracture, the iteration on the linearization is not necessary and does not affect the solution. Thus, we only compute one linearization \auth{step}. In these examples, we exclusively use the extrapolation.

In \RA{three} other examples, where we face growing fractures, the time-lagging phenomenon due to the linearization can be observed very well. In these \auth{tests}, we use the iteration on the linearization (ItL) until a certain tolerance is reached and compare the Newton iterations with and without ItL for both the iteration on the extrapolation (ItE) and the iteration via the previous incremental step solution (ItOTS). We also test the active set modifications from Section \ref{subsec:mod_newton_active_set}, but we did not observe significant differences between \textbf{Case 1} and \textbf{Case 2}. We assume, that the modification predominantly affects the number of active set iterations on very fine meshes (see Section \ref{subsec:sneddon_2d} and Section \ref{subsec:sneddon_3d}). Thus, all results shown in Section \ref{subsec:asym_tp} \auth{to} \ref{subsec:sens} are obtained from \textbf{Case 2}.

\subsection{Sneddon 2d}\label{subsec:sneddon_2d}
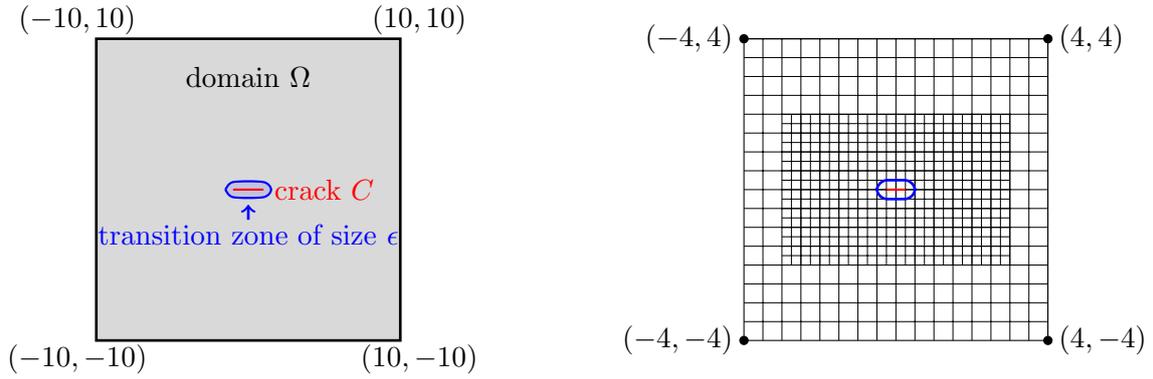
\begin{figure}[htbp!]
\begin{minipage}{0.47\textwidth}
\centering
\begin{tikzpicture}[scale = 1.0]
\node  at (-0.25,4.25) {$(-10,10)$};
\node  at (-0.25,-0.25) {$(-10,-10)$};
\node at (4.25,-0.25) {$(10,-10)$};
\node at (4.25,4.25) {$(10,10)$};
\draw[fill=gray!30] (0,0)  -- (0,4) -- (4,4)  -- (4,0) -- cycle;
\node at (2.0,3.5) {domain $\Omega$};
\draw [fill,opacity=0.1, blue] plot [smooth] coordinates { (1.7,2) (1.8,1.9) (2.2,1.9) (2.3,2) (2.2,2.1) (1.8,2.1) (1.7,2)};
\draw [thick, blue] plot [smooth] coordinates { (1.7,2) (1.8,1.9) (2.2,1.9) (2.3,2) (2.2,2.1) (1.8,2.1) (1.7,2)};
\draw[thick, draw=red] (1.8,2)--(2.2,2);
\node[red] at (3,2) {crack $C$};
\draw[->,blue] (2,1.6)--(2,1.8);
\node[blue] at (2,1.4) {transition zone of size $\epsilon$};
\end{tikzpicture}
    \end{minipage}
    \hspace{0.2cm}
    \begin{minipage}{0.47\textwidth}
    \centering
      \begin{tikzpicture}[scale = 2.5]
        \draw[black, thin] (-0.8,-0.8) -- (-0.8,0.8);
        \draw[black, thin] (-0.8,0.8) -- (0.8,0.8);
        \draw[black, thin] (0.8,0.8) -- (0.8,-0.8);
        \draw[black, thin] (0.8,-0.8) -- (-0.8,-0.8);
        \draw[step=0.1,black,ultra thin] (-0.8,-0.8) grid (0.8,0.8);
        \draw[step=0.05,black,ultra thin] (-0.6,-0.4) grid (0.6,0.4);
        \filldraw[black] (-0.8,-0.8) circle (0.5pt) node[anchor=east] {$(-4,-4)$};
        \filldraw[black] (-0.8,0.8) circle (0.5pt) node[anchor=east] {$(-4,4)$};
        \filldraw[black] (0.8,-0.8) circle (0.5pt) node[anchor=west] {$(4,-4)$};
        \filldraw[black] (0.8,0.8) circle (0.5pt) node[anchor=west] {$(4,4)$};
        \draw[thick, red] (-0.05,0) -- (0.05,0);
        \draw[blue, rounded corners=4pt] (-0.1,-0.05) rectangle ++(0.2,0.1);
    \end{tikzpicture}      
    \end{minipage}
    \caption{Left: geometry of the two dimensional Sneddon test. Right: zoom-in to the pre-refined crack zone in $[-4,4]\times[-4,4]$ with two global refinement steps and one local refinement step (geometrically pre-refined).}
    \label{geo_sneddon}
\end{figure}

We consider a stationary benchmark test~\cite{schroder2021selection}, where a constant pressure is applied in the inner of a pre-existing crack in the middle of a domain, and only the crack width varies. This test setup is motivated by Sneddon~\cite{sneddon1946distribution}, and Sneddon and Lowegrub~\cite{SneddLow69}. We restrict ourselves to a one dimensional fracture $C$ on a two dimensional domain $\Omega = (-10,10)^2$ as depicted in Figure~\ref{geo_sneddon} on the left. The fracture is centered horizontally within $\Omega$ and has a constant half crack length $l_0 = 0.25$ and varying width. Precisely, the crack width corresponds to $2h$, where $h$ is the minimal element diameter of the mesh. The mesh is pre-refined geometrically in the crack zone, as depicted exemplarily for one adaptive refinement step in Figure~\ref{geo_sneddon} on the right, where the crack zone is resolved with the smallest mesh size. The driving force is given by a constant pressure $\rho = 10^{-3}\,\mathrm{Pa}$ in the inner crack. The parameter setting is given in Table~\ref{table_param}.

The spatial discretization parameter, i.e. the minimal element diameter, is set as
\[
h=0.022, \, 0.011, \, 0.0055, \, 0.0027, \, 0.0013.
\]
The quantity of interest, called total crack volume ($\operatorname{TCV}$), can be computed numerically via
\begin{align*}
\operatorname{TCV} = \int_\Omega u(x,y) \nabla \varphi (x,y)\, d(x,y).
\end{align*}
The analytical solution~\cite{SneddLow69} is given by
\begin{align*}
\operatorname{TCV}_{\operatorname{ref}} = \frac{2\pi \rho l_0^2}{E'}.
\end{align*}

\begin{table}[htbp!]
\renewcommand{\arraystretch}{1.2}
\scriptsize
\centering
\begin{tabular}{|l|l|c|}\hline
\multicolumn{1}{|c|}{Parameter} & \multicolumn{1}{c|}{Definition} &  \multicolumn{1}{c|}{Value} \\ \hline
$\Omega$ & Domain & $(-10,10)^2$  \\ 
$ h $ & Diagonal cell diameter & test-dependent\\
$ l_0 $ & Half crack length & $0.25$ \\ 
$ G_C $ & Material toughness & $1.0$ \\ 
$ E $ & Young's modulus & $1.0$ \\ 
$\mu$ & Lamé parameter &  $0.42$\\ 
$\lambda$ & Lamé parameter & $0.28$ \\ 
$ \nu $ & Poisson's ratio & $0.2$ \\ 
$ p $ & Applied pressure & $10^{-3}$ \\ 
$ \varepsilon $ & Bandwidth of the initial crack & $2h$ \\
$ \kappa $ & Regularization parameter & $10^{-10}$\\ 
& \RA{Number of global refinements} & \RA{$2$}\\
& \RA{Number of local refinements} & \RA{$5$, $6$, $7$, $8$, $9$}\\
\RA{$\operatorname{TOL}_N$}& \RA{Tolerance outer Newton solver} & \RA{$10^{-7}$}  \\
& \RA{Tolerance inner linear solver} & \RA{$\lVert \Tilde{R}(U_h^{n,k})\rVert_2 10^{-8}$}  \\\hline 
\end{tabular}
\caption{The setting of the material and numerical parameters for the Sneddon 2d test.}
\label{table_param}
\end{table}

The Figures \ref{plot_as_iter_0.088} - \ref{plot_as_iter_0.011} visualize the average active-set iterations required per time-step for \textbf{Case 1}-\textbf{Case 4} and $h=0.022,0.011,0.0055,0.0027$. We observe, that with a smaller $h$, the number of iterations can become comparatively large in \textbf{Case 1} in comparison to the other cases. We also ran a test with a minimal element diameter $h=0.0013$ which leads to around 40 million degrees of freedom. There, the Newton active set method did not achieve convergence within $500$ Newton active set iterations with \textbf{Case 1} in the $2$nd incremental step. In contrast, with the modified constant $c$ from \textbf{Case 2}, the Newton method terminates within 
$46,\, 26,\, 18, \, 7$ and $4$ iterations to converge at each incremental step. With \textbf{Case 3} $10-11$ Newton active set iterations are needed in each timestep on each refinement level. With \textbf{Case 4} only $1$ iteration is needed. The latter shows that the Newton method only needs $1$ iteration to converge, which underlines the assumption that the active set stopping criterion is indeed the reason for the slow convergence in \textbf{Case 1}. 
\RA{The speed-up coming along with \textbf{Case 3}, \textbf{Case 4} and especially \textbf{Case 2} has a noticeable impact in the computation time, which is summarized for all settings in Table \ref{table_wall_time_sneddon_2d}. The speed-up of \textbf{Case 3} and \textbf{Case 4} must be treated with caution since we do not iterate until full convergence of the active set in these cases, which may lead to reduced accuracy.}

\begin{table}[htbp!]
\renewcommand{\arraystretch}{1.2}
\scriptsize
\centering
\RA{
\begin{tabular}{|l|c|c|c|c|c|}\hline
& \multicolumn{5}{c|}{Total wallclock time $|$ Sneddon 2d}  \\ \hline
\multicolumn{1}{|c|}{Case} & \multicolumn{1}{c|}{$h=0.022$ $|$ 4 cores} & \multicolumn{1}{c|}{$h=0.011$ $|$ 4 cores} & \multicolumn{1}{c|}{$h=0.0055$ $|$ 16 cores} & \multicolumn{1}{c|}{$h=0.0027$ $|$ 32 cores} & \multicolumn{1}{c|}{$h=0.0013$ $|$ 64 cores}\\\hline
Case 1 & $587.983\si{s}$ & $2243.164\si{s}$ & $2275.481\si{s}$ & $11682.962\si{s}$ & --\\\hline
Case 2 & $161.765\si{s}$ & $976.391\si{s}$ & $784.019\si{s}$ & $2590.394\si{s}$ & $6562.554\si{s}$\\\hline
Case 3 & $117.775\si{s}$ & $458.664\si{s}$ & $572.438\si{s}$ & $1228.499\si{s}$ & $4132.350\si{s}$\\\hline
Case 4 & $22.833\si{s}$  & $77.788\si{s}$ & $87.405\si{s}$ & $208.310\si{s}$ & $940.686\si{s}$\\\hline
\end{tabular}
}
\caption{\RA{Total wallclock time of the Sneddon $2$d test for different refinement levels on different numbers of cores.}}
\label{table_wall_time_sneddon_2d}
\end{table}

As it can be observed in Table \ref{table_tcv_err_2d}, the error in the TCV is not affected, even though we completely ignore the active set stopping criterion (\textbf{Case 4}). But in this case, we may obtain a non-smooth phase-field solution. This phenomenon is depicted in Figure \ref{fig:nonsmooth_pf}. The slight increase of the TCV error for $h<0.0055$ can be attributed to the fact that the analytical solution is based on an infinite domain. Thus, the solution does not converge to the exact solution for $h\rightarrow 0$ since the domain-error due to the finite domain will become dominant if $h$ is small enough, which is the case for $h<0.0055$, as we assume. This phenomenon was further investigated in \cite{HeiWi18_pamm}.

\begin{figure}[htbp!]
\begin{minipage}[b]{0.49\textwidth}
\scriptsize
\centering
\begin{tikzpicture}[scale = 1.0]
\begin{axis}[
    xlabel = time, 
    ylabel = $\#$ Active set iterations,
    legend style = {at={(0.75,1.3)}},
    mark size=0.3pt,
    grid =major,
    y post scale = 0.9,
    xtick={0,1,2,3,4}
    ]
        \addplot[line width = 6pt,color=blue, mark=x]
        table[col sep=space] {sneddon_as_iter_case_1_0.022.txt};
        \addlegendentry{\textbf{Case 1} $|$ $h=0.022$}
        \addplot[line width = 6pt,color=red, mark=x, dashed]
        table[col sep=space] {sneddon_as_iter_case_2_0.022.txt};
        \addlegendentry{\textbf{Case 2} $|$ $h=0.022$}
        \addplot[line width = 6pt,color=green, mark=x, densely dotted]
        table[col sep=space] {sneddon_as_iter_case_3_0.022.txt};
        \addlegendentry{\textbf{Case 3} $|$ $h=0.022$}
        \addplot[line width = 6pt,color=brown, mark=x, densely dashed]
        table[col sep=space] {sneddon_as_iter_case_4_0.022.txt};
        \addlegendentry{\textbf{Case 4} $|$ $h=0.022$}
\end{axis}
\end{tikzpicture}
\caption{Number of active set iterations for \textbf{Case 1}-\textbf{Case 4} and $h=0.022$ in the Sneddon 2d test.}
\label{plot_as_iter_0.088}

\end{minipage}
\hfill
\begin{minipage}[b]{0.49\textwidth}
\scriptsize
\centering
\begin{tikzpicture}[scale = 1.0]
\begin{axis}[
    xlabel = time, 
    ylabel = $\#$ Active set iterations,
    legend style = {at={(0.75,1.3)}},
    mark size=0.3pt,
    grid =major,
    y post scale = 0.9,
    xtick={0,1,2,3,4}
    ]
        \addplot[line width = 6pt,color=blue, mark=x]
        table[col sep=space] {sneddon_as_iter_case_1_0.011.txt};
        \addlegendentry{\textbf{Case 1} $|$ $h=0.011$}
        \addplot[line width = 6pt,color=red, mark=x, dashed]
        table[col sep=space] {sneddon_as_iter_case_2_0.011.txt};
        \addlegendentry{\textbf{Case 2} $|$ $h=0.011$}
        \addplot[line width = 6pt,color=green, mark=x, densely dotted]
        table[col sep=space] {sneddon_as_iter_case_3_0.011.txt};
        \addlegendentry{\textbf{Case 3} $|$ $h=0.011$}
        \addplot[line width = 6pt,color=brown, mark=x, densely dashed]
        table[col sep=space] {sneddon_as_iter_case_4_0.011.txt};
        \addlegendentry{\textbf{Case 4} $|$ $h=0.011$}
\end{axis}
\end{tikzpicture}
\caption{Number of active set iterations for \textbf{Case 1} to \textbf{Case 4} and $h=0.011$ in the Sneddon 2d test.}
\label{plot_as_iter_0.044}
\end{minipage}
\end{figure}

\begin{figure}[htbp!]
\begin{minipage}[b]{0.49\textwidth}
\scriptsize
\centering
\begin{tikzpicture}[scale = 1.0]
\begin{axis}[
    xlabel = time, 
    ylabel = $\#$ Active set iterations,
    legend style = {at={(0.75,1.3)}},
    mark size=0.3pt,
    grid =major,
    y post scale = 0.9,
    xtick={0,1,2,3,4}
    ]
        \addplot[line width = 6pt,color=blue, mark=x]
        table[col sep=space] {sneddon_as_iter_case_1_0.0055.txt};
        \addlegendentry{\textbf{Case 1} $|$ $h=0.0055$}
        \addplot[line width = 6pt,color=red, mark=x, dashed]
        table[col sep=space] {sneddon_as_iter_case_2_0.0055.txt};
        \addlegendentry{\textbf{Case 2} $|$ $h=0.0055$}
        \addplot[line width = 6pt,color=green, mark=x, densely dotted]
        table[col sep=space] {sneddon_as_iter_case_3_0.0055.txt};
        \addlegendentry{\textbf{Case 3} $|$ $h=0.0055$}
        \addplot[line width = 6pt,color=brown, mark=x, densely dashed]
        table[col sep=space] {sneddon_as_iter_case_4_0.0055.txt};
        \addlegendentry{\textbf{Case 4} $|$ $h=0.0055$}
\end{axis}
\end{tikzpicture}
\caption{Number of active set iterations for \textbf{Case 1}-\textbf{Case 4} and $h=0.0055$ in the Sneddon 2d test.}
\label{plot_as_iter_0.022}
\end{minipage}
\hfill
\begin{minipage}[b]{0.49\textwidth}
\scriptsize
\centering
\begin{tikzpicture}[scale = 1.0]
\begin{axis}[
    xlabel = time, 
    ylabel = $\#$ Active set iterations,
    legend style = {at={(0.75,1.3)}},
    mark size=0.3pt,
    grid =major,
    y post scale = 0.9,
    xtick={0,1,2,3,4}
    ]
        \addplot[line width = 6pt,color=blue, mark=x]
        table[col sep=space] {sneddon_as_iter_case_1_0.0027.txt};
        \addlegendentry{\textbf{Case 1} $|$ $h=0.0027$}
        \addplot[line width = 6pt,color=red, mark=x, dashed]
        table[col sep=space] {sneddon_as_iter_case_2_0.0027.txt};
        \addlegendentry{\textbf{Case 2} $|$ $h=0.0027$}
        \addplot[line width = 6pt,color=green, mark=x, densely dotted]
        table[col sep=space] {sneddon_as_iter_case_3_0.0027.txt};
        \addlegendentry{\textbf{Case 3} $|$ $h=0.0027$}
        \addplot[line width = 6pt,color=brown, mark=x, densely dashed]
        table[col sep=space] {sneddon_as_iter_case_4_0.0027.txt};
        \addlegendentry{\textbf{Case 4} $|$ $h=0.0027$}
\end{axis}
\end{tikzpicture}
\caption{Number of active set iterations for \textbf{Case 1}-\textbf{Case 4} and $h=0.0027$ in the Sneddon 2d test.}
\label{plot_as_iter_0.011}
\end{minipage}
\end{figure}

\begin{table}[htbp!]
\renewcommand{\arraystretch}{1.2}
\scriptsize
\centering
\begin{tabular}{|l|c|c|c|c|c|}\hline
& \multicolumn{5}{c|}{TCV error $|$ Sneddon 2d}  \\ \hline
\multicolumn{1}{|c|}{Case} & \multicolumn{1}{c|}{$h=0.022$} & \multicolumn{1}{c|}{$h=0.011$} & \multicolumn{1}{c|}{$h=0.0055$} & \multicolumn{1}{c|}{$h=0.0027$} & \multicolumn{1}{c|}{$h=0.0013$}\\\hline
Case 1 & 0.000173649  & 3.72255e-05 & 3.03573e-05 & 6.45367e-05 & --\\\hline
Case 2 & 0.000173649  & 3.72255e-05 & 3.03573e-05 & 6.45367e-05 & 8.19169e-05\\\hline
Case 3 & 0.000173649  & 3.72255e-05 & 3.03573e-05 & 6.45367e-05 & 8.19169e-05\\\hline
Case 4 & 0.000173649  & 5.33208e-05 & 3.03573e-05 & 6.14307e-05 & 8.19169e-05\\\hline
\end{tabular}
\caption{Error in the TCV for different mesh size parameters and \textbf{Case 1} - {Case 4}.}
\label{table_tcv_err_2d}
\end{table}
\begin{figure}
    \centering
    \includegraphics[scale = 0.24]{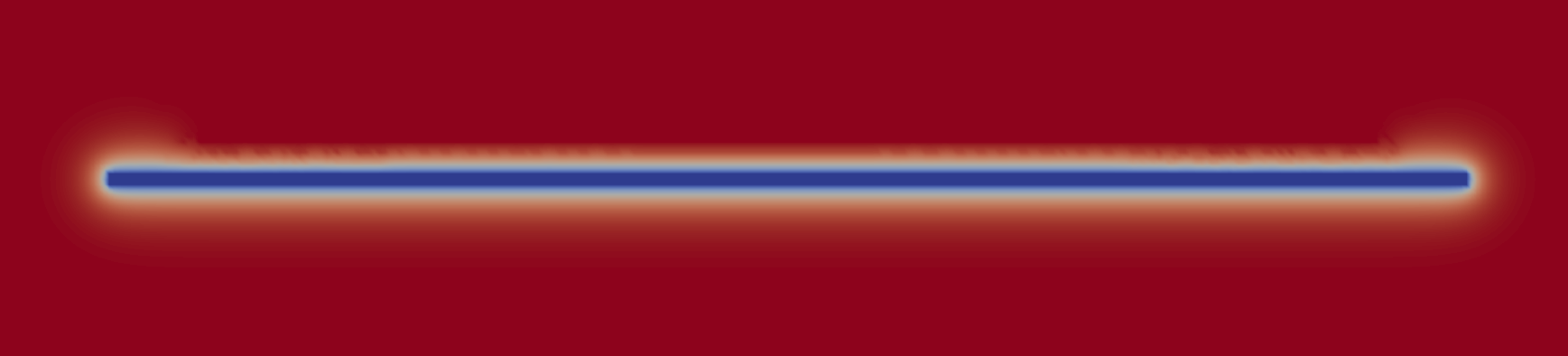}
    \caption{Visualization of the nonsmooth phase-field solution of the two dimensional Sneddon test with $h=0.011$ and \textbf{Case 4}. \RA{Red represents the fully intact part of the domain, blue the fully broken part and white stands for the transition zone.}}
    \label{fig:nonsmooth_pf}
\end{figure}

Lastly, a brief parallel study is conducted. The framework is parallelized with MPI. 
Since it is based on the \texttt{pfm-cracks} code \cite{heister2020pfm-cracks}, extensive scalability analyses can be found in \cite{HeiWi18_pamm}. \RA{We compare the CPU time for one representative incremental step for $h=0.022$ ($168\, 609$ degrees of freedom) computed on $1$ core and $16$ cores. On $1$ core, the CPU time of one incremental step is $2386\si{s}$ (approx. $18\si{h}$) and on $16$ cores, the CPU time is $142\si{s}$.}

\subsection{Sneddon 3d}\label{subsec:sneddon_3d}
The Sneddon 3d test \cite{SneddLow69}[Section 3.3] is the three dimensional equivalent of the Sneddon 2d test. We have a three dimensional cubic domain $\Omega$ and a two dimensional penny-shaped fracture with a pressure acting in its center. As in Section \ref{subsec:sneddon_2d}, we compare the number of Newton iterations for \textbf{Case 1} - \textbf{Case 4}. The results for different mesh size parameters $h = 1.732,\, 0.866,\, 0.433,\, 0.216$ are depicted in Figure \ref{plot_sneddon_3d_as_iter_1.732} - Figure \ref{plot_sneddon_3d_as_iter_0.216}. On coarser grids, \textbf{Case 1} and \textbf{Case 2} perform very similarly, whereas on finer meshes, the benefit of \textbf{Case 2} is seen especially in the fourth incremental step. As assumed, \textbf{Case 4} yields the lowest number of active set iterations and \textbf{Case 3} yields a constant number of iterations of around $10-12$ per incremental step.
\begin{table}[htbp!]
\renewcommand{\arraystretch}{1.2}
\scriptsize
\centering
\begin{tabular}{|l|l|c|}\hline
\multicolumn{1}{|c|}{Parameter} & \multicolumn{1}{c|}{Definition} &  \multicolumn{1}{c|}{Value} \\ \hline
$\Omega$ & Domain & $(-10,10)^3$  \\ 
$ h $ & Diagonal cell diameter & test-dependent\\
$ r_0 $ & crack radius & $1.0$ \\ 
$ G_C $ & Material toughness & $1.0$ \\ 
$ E $ & Young's modulus & $1.0$ \\ 
$\mu$ & Lamé parameter &  $0.42$\\ 
$\lambda$ & Lamé parameter & $0.28$ \\ 
$ \nu $ & Poisson's ratio & $0.2$ \\ 
$ p $ & Applied pressure & $10^{-3}$ \\ 
$ \varepsilon $ & Bandwidth of the initial crack & $2h$ \\
$ \kappa $ & Regularization parameter & $10^{-10}h$\\ 
& \RA{Number of global refinements} & \RA{$1$}\\
& \RA{Number of local refinements} & \RA{$0$, $1$, $2$, $3$}\\
\RA{$\operatorname{TOL}_N$}& \RA{Tolerance outer Newton solver} & \RA{$10^{-7}$}  \\ 
& \RA{Tolerance inner linear solver} & \RA{$\lVert \Tilde{R}(U_h^{n,k})\rVert_2 10^{-8}$}  \\\hline 
\end{tabular}
\caption{The setting of the material and numerical parameters used for the Sneddon 3d test.}
\label{table_param_3d}
\end{table}

\begin{figure}[htbp!]
\begin{minipage}[b]{0.49\textwidth}
\scriptsize
\centering
\begin{tikzpicture}[scale = 1.0]
\begin{axis}[
    xlabel = time, 
    xmax = 5,
    ylabel = $\#$ Active set iterations,
    legend style = {at={(0.75,1.3)}},
    mark size=0.3pt,
    grid =major,
    y post scale = 0.9,
    xtick={0.0,1.0,2.0,3.0,4.0,5.0}
    ]
        \addplot[line width = 6pt,color=blue, mark=x]
        table[col sep=space] {sneddon_3d_as_iter_case_1_1.732.txt};
        \addlegendentry{\textbf{Case 1} $|$ $h=1.732$}
        \addplot[line width = 6pt,color=red, mark=x, dashed]
        table[col sep=space] {sneddon_3d_as_iter_case_2_1.732.txt};
        \addlegendentry{\textbf{Case 2} $|$ $h=1.732$}
        \addplot[line width = 6pt,color=green, mark=x, densely dotted]
        table[col sep=space] {sneddon_3d_as_iter_case_3_1.732.txt};
        \addlegendentry{\textbf{Case 3} $|$ $h=1.732$}
        \addplot[line width = 6pt,color=brown, mark=x, densely dashed]
        table[col sep=space] {sneddon_3d_as_iter_case_4_1.732.txt};
        \addlegendentry{\textbf{Case 4} $|$ $h=1.732$}
\end{axis}
\end{tikzpicture}
\caption{Number of active set iterations for \textbf{Case 1}-\textbf{Case 4} and $h=1.732$ for the Sneddon 3d test.}
\label{plot_sneddon_3d_as_iter_1.732}

\end{minipage}
\hfill
\begin{minipage}[b]{0.49\textwidth}
\scriptsize
\centering
\begin{tikzpicture}[scale = 1.0]
\begin{axis}[
    xlabel = time, 
    xmax = 5,
    ylabel = $\#$ Active set iterations,
    legend style = {at={(0.75,1.3)}},
    mark size=0.3pt,
    grid =major,
    y post scale = 0.9,
    xtick={0.0,1.0,2.0,3.0, 4.0,5.0}
    ]
        \addplot[line width = 6pt,color=blue, mark=x]
        table[col sep=space] {sneddon_3d_as_iter_case_1_0.866.txt};
        \addlegendentry{\textbf{Case 1} $|$ $h=0.866$}
        \addplot[line width = 6pt,color=red, mark=x, dashed]
        table[col sep=space] {sneddon_3d_as_iter_case_2_0.866.txt};
        \addlegendentry{\textbf{Case 2} $|$ $h=0.866$}
        \addplot[line width = 6pt,color=green, mark=x, densely dotted]
        table[col sep=space] {sneddon_3d_as_iter_case_3_0.866.txt};
        \addlegendentry{\textbf{Case 3} $|$ $h=0.866$}
        \addplot[line width = 6pt,color=brown, mark=x, densely dashed]
        table[col sep=space] {sneddon_3d_as_iter_case_4_0.866.txt};
        \addlegendentry{\textbf{Case 4} $|$ $h=0.866$}
\end{axis}
\end{tikzpicture}
\caption{Number of active set iterations for \textbf{Case 1}-\textbf{Case 4} and $h=0.866$ for the Sneddon 3d test.}
\label{plot_sneddon_3d_as_iter_0.866}
\end{minipage}
\end{figure}

\begin{figure}[htbp!]
\begin{minipage}[b]{0.49\textwidth}
\scriptsize
\centering
\begin{tikzpicture}[scale = 1.0]
\begin{axis}[
    xlabel = time, 
    xmax = 5,
    ylabel = $\#$ Active set iterations,
    legend style = {at={(0.75,1.3)}},
    mark size=0.3pt,
    grid =major,
    y post scale = 0.9,
    xtick={0.0,1.0,2.0,3.0,4.0,5.0}
    ]
        \addplot[line width = 6pt,color=blue, mark=x]
        table[col sep=space] {sneddon_3d_as_iter_case_1_0.433.txt};
        \addlegendentry{\textbf{Case 1} $|$ $h=0.433$}
        \addplot[line width = 6pt,color=red, mark=x, dashed]
        table[col sep=space] {sneddon_3d_as_iter_case_2_0.433.txt};
        \addlegendentry{\textbf{Case 2} $|$ $h=0.433$}
        \addplot[line width = 6pt,color=green, mark=x, densely dotted]
        table[col sep=space] {sneddon_3d_as_iter_case_3_0.433.txt};
        \addlegendentry{\textbf{Case 3} $|$ $h=0.433$}
        \addplot[line width = 6pt,color=brown, mark=x, densely dashed]
        table[col sep=space] {sneddon_3d_as_iter_case_4_0.433.txt};
        \addlegendentry{\textbf{Case 4} $|$ $h=0.433$}
\end{axis}
\end{tikzpicture}
\caption{Number of active set iterations for \textbf{Case 1}-\textbf{Case 4} and $h=0.433$ for the Sneddon 3d test.}
\label{plot_sneddon_3d_as_iter_0.433}

\end{minipage}
\hfill
\begin{minipage}[b]{0.49\textwidth}
\scriptsize
\centering
\begin{tikzpicture}[scale = 1.0]
\begin{axis}[
    xlabel = time, 
    xmax = 5,
    ylabel = $\#$ Active set iterations,
    legend style = {at={(0.75,1.3)}},
    mark size=0.3pt,
    grid =major,
    y post scale = 0.9,
    xtick={0.0,1.0,2.0,3.0, 4.0,5.0}
    ]
        \addplot[line width = 6pt,color=blue, mark=x]
        table[col sep=space] {sneddon_3d_as_iter_case_1_0.216.txt};
        \addlegendentry{\textbf{Case 1} $|$ $h=0.216$}
        \addplot[line width = 6pt,color=red, mark=x, dashed]
        table[col sep=space] {sneddon_3d_as_iter_case_2_0.216.txt};
        \addlegendentry{\textbf{Case 2} $|$ $h=0.216$}
        \addplot[line width = 6pt,color=green, mark=x, densely dotted]
        table[col sep=space] {sneddon_3d_as_iter_case_3_0.216.txt};
        \addlegendentry{\textbf{Case 3} $|$ $h=0.216$}
        \addplot[line width = 6pt,color=brown, mark=x, densely dashed]
        table[col sep=space] {sneddon_3d_as_iter_case_4_0.216.txt};
        \addlegendentry{\textbf{Case 4} $|$ $h=0.216$}
\end{axis}
\end{tikzpicture}
\caption{Number of active set iterations for \textbf{Case 1}-\textbf{Case 4} and $h=0.216$ for the Sneddon 3d test.}
\label{plot_sneddon_3d_as_iter_0.216}
\end{minipage}
\end{figure}
\begin{table}[htbp!]
\renewcommand{\arraystretch}{1.2}
\scriptsize
\centering
\begin{tabular}{|l|c|c|c|c|}\hline
& \multicolumn{4}{c|}{TCV error $|$ Sneddon 3d}  \\ \hline
\multicolumn{1}{|c|}{Case} & \multicolumn{1}{c|}{$h=1.732$} & \multicolumn{1}{c|}{$h=0.866$} & \multicolumn{1}{c|}{$h=0.433$} & \multicolumn{1}{c|}{$h=0.216$}\\\hline
Case 1 & 0.061265 & 0.0220069 & 0.00897784 & 0.00352026\\\hline
Case 2 & 0.061265 & 0.0220069 & 0.00897784 & 0.00352026\\\hline
Case 3 & 0.061265 & 0.0220069 & 0.00897784 & 0.00352026\\\hline
Case 4 & 0.061265 & 0.0220069 & 0.00897606 & 0.00352025\\\hline
\end{tabular}
\caption{Error in the TCV of the $3$ dimensional Sneddon test for different mesh size parameters and \textbf{Case 1} - \textbf{Case 4}.}
\label{table_tcv_err_3d}
\end{table}

We mention since the crack does not grow in Section 
\ref{subsec:sneddon_2d} and Section \ref{subsec:sneddon_3d} that there is 
numerically no difference whether ItL is utilized or not. For this reason,
we concentrate on the active set performance only. \RA{However, as in the two dimensional case, the reduction of the number of Newton iterations with \textbf{Case 2}, \textbf{Case 3} and \textbf{Case 4}
leads to a significant speed-up, which can be observed by comparing the total wallclock times and is summarized in Table \ref{table_wall_time_sneddon_3d}.}
\begin{table}[htbp!]
\renewcommand{\arraystretch}{1.2}
\scriptsize
\centering
\RA{
\begin{tabular}{|l|c|c|c|c|c|}\hline
& \multicolumn{4}{c|}{Total wallclock time $|$ Sneddon 3d}  \\ \hline
\multicolumn{1}{|c|}{Case} & \multicolumn{1}{c|}{$h=1.732$ $|$ 32 cores} & \multicolumn{1}{c|}{$h=0.866$ $|$ 48 cores} & \multicolumn{1}{c|}{$h=0.433$ $|$ 128 cores} & \multicolumn{1}{c|}{$h=0.216$ $|$ 128 cores}\\\hline
Case 1 & $19.669\si{s}$ & $24.089\si{s}$ & $267.847\si{s}$ & $888.671\si{s}$ \\\hline
Case 2 & $19.981\si{s}$ & $21.836\si{s}$ & $191.004\si{s}$ & $637.236\si{s}$\\\hline
Case 3 & $36.249\si{s}$ & $37.906\si{s}$ & $322.549\si{s}$ & $590.830\si{s}$ \\\hline
Case 4 & $9.949\si{s}$  & $10.580\si{s}$ & $87.367\si{s}$ & $162.499\si{s}$ \\\hline
\end{tabular}
}
\caption{\RA{Total wallclock time of the Sneddon 3d test for different refinement levels on different numbers of cores.}}
\label{table_wall_time_sneddon_3d}
\end{table}

Finally, in this example, we \auth{test} the parallel performance with $h=0.433$ \RA{by comparing the CPU time for one incremental step}. 
On $1$ core, \RA{$8889\si{s}$ (approx.\ $2.4\si{h}$)}
were needed 
whereas the 
\RA{CPU time}
on \RA{$16$} cores is \RA{$542\si{s}$} (approx. \RA{$9$} minutes).

\subsection{Asymmetric three-point bending test}\label{subsec:asym_tp}
\begin{figure}[htbp!]
\centering
\begin{tikzpicture}[xscale=0.6,yscale=0.6]
\draw[fill=gray!30] (0,0) -- (0,8) -- (20,8) -- (20,0) -- (0,0);
\draw[blue] (20,8) -- (0,8);
\draw[red] (4,0) -- (4,1.0);
\node[red] at (3.1,1.0) {slit};
\draw[black] (6.0,2.75) circle (0.15);
\draw[black] (6.0,4.75) circle (0.15);
\draw[black] (6.0,6.75) circle (0.15);
\draw[->,blue] (10,10.5) -- (10,8.5);
\draw[->,blue] (6,10.5) -- (6,9);
\draw[->,blue] (14,10.5) -- (14,9);
\draw[->,blue] (2,10.5) -- (2,9.5);
\draw[->,blue] (18,10.5) -- (18,9.5);
\node at (9.5,10.5) {$u_y$};
\node at (20.15,-0.5) {$x$};
\node at (-0.5,8.15) {$y$};

\draw[<->] (20.5,0) -- (20.5,8);
\node[anchor=west] at (20.5,4.0) {$8\,\mathrm{mm}$};

\draw[<->] (7,0) -- (7,2.6);
\node[anchor=west] at (7,1.3) {$2.6\,\mathrm{mm}$};

\draw[<->] (7,2.9) -- (7,4.6);
\node[anchor=west] at (7,3.75) {$1.7\,\mathrm{mm}$};

\draw[<->] (0,3.75) -- (5.85,3.75);
\node[anchor=south] at (2.925,3.75) {$5.85\,\mathrm{mm}$};

\draw[<->] (7,4.9) -- (7,6.6);
\node[anchor=west] at (7,5.75) {$1.7\,\mathrm{mm}$};

\draw[<->] (7,6.9) -- (7,8);
\node[anchor=west] at (7,7.45) {$1.1\,\mathrm{mm}$};

\draw[<->] (4.5,0) -- (4.5,1);
\node[anchor=west] at (4.5,0.5) {$1\,\mathrm{mm}$};

\draw[black] (19,0) -- (17.5,-1.5);
\draw[black] (17.5,-1.5) -- (20.5,-1.5);
\draw[black] (20.5,-1.5) -- (19,0);
\draw[black] (17.5,-1.8) circle (0.3);
\draw[black] (20.5,-1.8) circle (0.3);
\draw[black] (18.1,-1.8) circle (0.3);
\draw[black] (19.9,-1.8) circle (0.3);
\draw[black] (18.7,-1.8) circle (0.3);
\draw[black] (19.3,-1.8) circle (0.3);

\draw[black] (1,0) -- (-0.5,-1.5);
\draw[black] (-0.5,-1.5) -- (2.5,-1.5);
\draw[black] (2.5,-1.5) -- (1,0);
\draw[black] (-0.5,-1.5) -- (-0.8,-1.7);
\draw[black] (-0.2,-1.5) -- (-0.5,-1.7);
\draw[black] (0.1,-1.5) -- (-0.2,-1.7);
\draw[black] (0.4,-1.5) -- (0.1,-1.7);
\draw[black] (0.7,-1.5) -- (0.4,-1.7);
\draw[black] (1.0,-1.5) -- (0.7,-1.7);
\draw[black] (1.3,-1.5) -- (1.0,-1.7);
\draw[black] (1.6,-1.5) -- (1.3,-1.7);
\draw[black] (1.9,-1.5) -- (1.6,-1.7);
\draw[black] (2.2,-1.5) -- (1.9,-1.7);
\draw[black] (2.5,-1.5) -- (2.2,-1.7);

\draw[<->] (0,-3) -- (1,-3);
\node[anchor=north] at (0.5,-3) {$1\,\mathrm{mm}$};

\draw[<->] (20,-3) -- (19,-3);
\node[anchor=north] at (19.5,-3) {$1\,\mathrm{mm}$};

\draw[<->] (0,-4.5) -- (20,-4.5);
\node[anchor=north] at (10,-4.5) {$20\,\mathrm{mm}$};
\end{tikzpicture}
 \caption{Visualization of the configuration of the asymmetric three-point bending test.}\label{fig:geo_astp}
\end{figure}
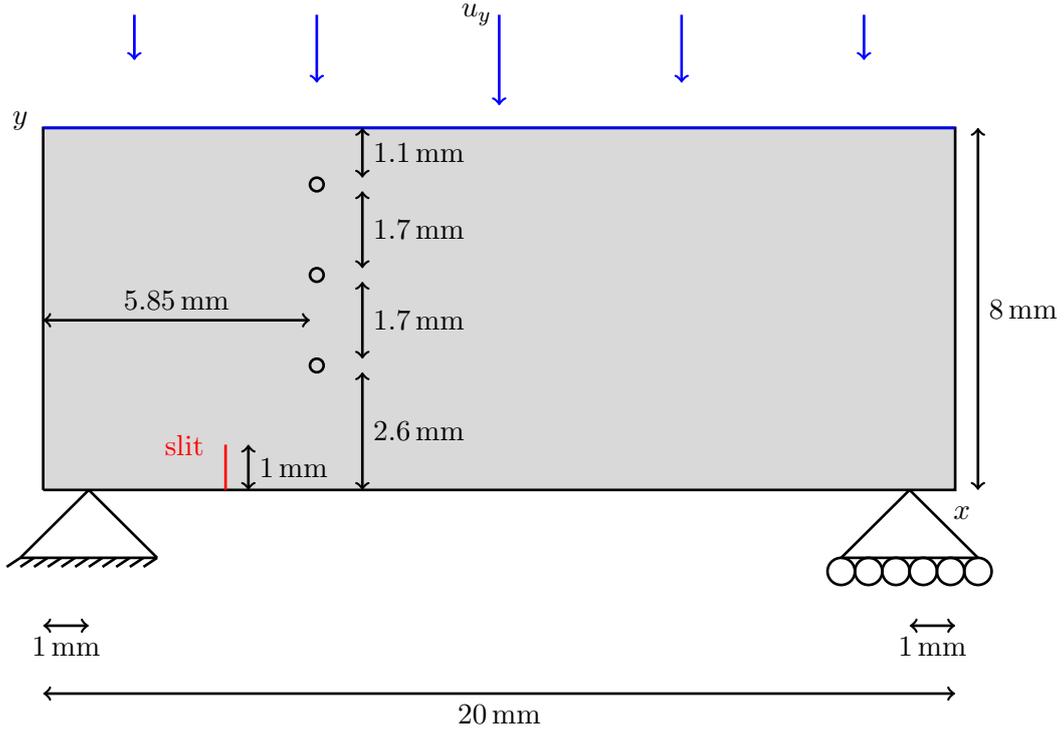

\begin{figure}[htbp!]
\begin{minipage}[b]{0.49\textwidth}
    \centering
    \includegraphics[scale = 0.20]{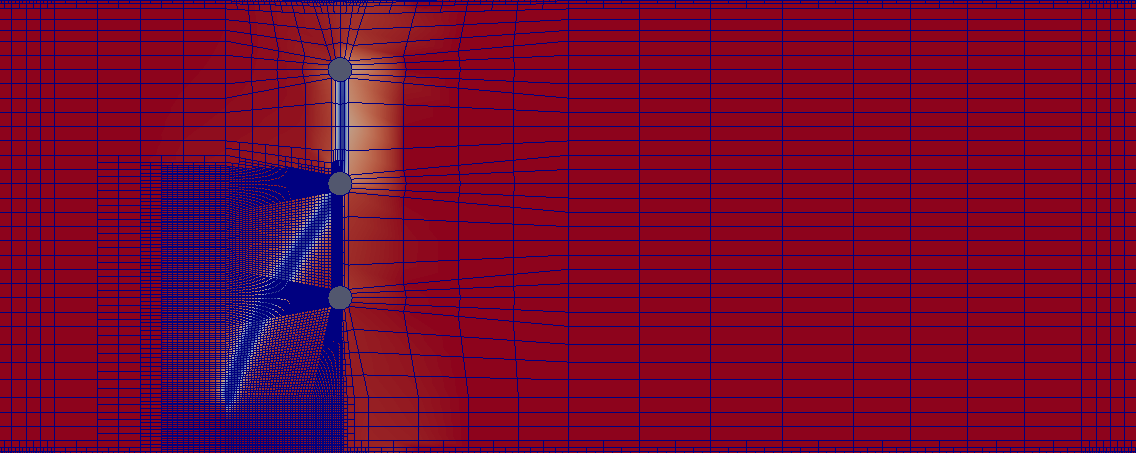}
    \caption{\RA{Visualization of the final configuration of the fracture asymmetric three-point bending test with $2$ global an $3$ local refinements.} \RA{Red represents the fully intact part of the domain, blue the fully broken part and white stands for the transition zone.}}
    \label{fig:astp_vis_final_conf_less_ref}
\end{minipage}
\hfill
\begin{minipage}[b]{0.49\textwidth}
    \centering
    \includegraphics[scale=0.20]{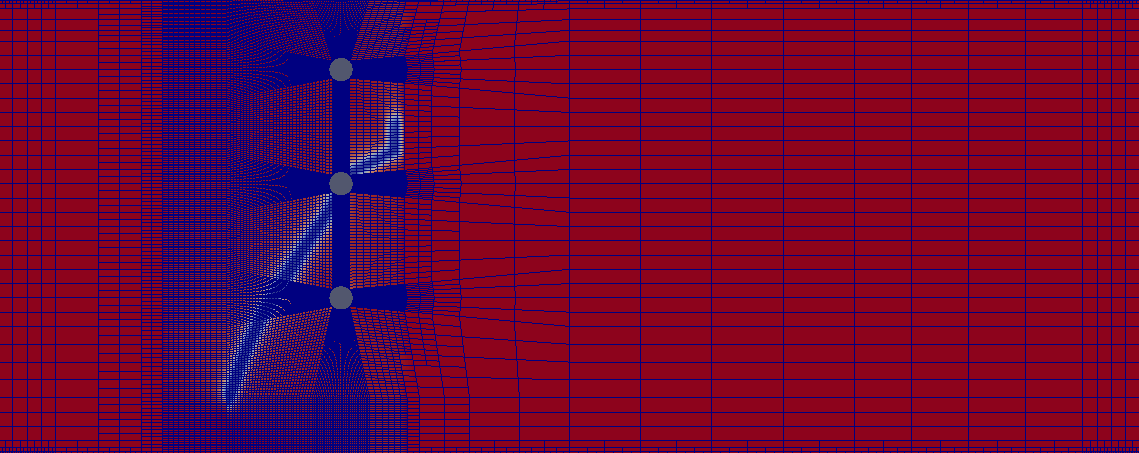}
    \caption{\RA{Visualization of the final configuration of the fracture asymmetric three-point bending test $2$ global an $3$ local refinements.} \RA{Red represents the fully intact part of the domain, blue the fully broken part and white stands for the transition zone.}}
    \label{fig:astp_vis_final_conf_more_ref}
\end{minipage}
\end{figure}

In this third example, we consider the asymmetric three-point bending test 
\cite{MieWelHof10a,MesBouKhon15,AmGeraLoren15}.
Here, the fracture grows and we examine the performance of the iteration on the linearization and its influence on the solution. The configuration is displayed in Figure \ref{fig:geo_astp}. 
We consider a \auth{two} dimensional domain with three holes. The boundary conditions are taken from 
\cite{BrWiBeNoRa20}. 
On $\partial \Omega_{\operatorname{top}}$, the upper boundary, we apply time-dependent non-homogeneous Dirichlet conditions in $y$-direction in a Gaussian bell curve fashion:
\begin{equation*}
    u_y  = -10\exp\left(-\frac{(x-10)^2}{100}\right)t\cdot 1\mathrm{mm/s}, \; x \in [0,10], \, t \in [0,T],
\end{equation*}
where $T>0$ is the maximum time. Furthermore, we fix the displacement in both directions at $(1,0)$ and in $y$-direction at $(19,0)$. Otherwise, the boundary conditions are defined to be traction-free (homogeneous Neumann conditions). The 
material and numerical
parameters for this test are \auth{provided} in Table \ref{table:params_astp}.
\begin{table}[htbp!]
\renewcommand{\arraystretch}{1.2}
\scriptsize
\centering
\begin{tabular}{|l|l|c|}\hline
\multicolumn{1}{|c|}{Parameter} & \multicolumn{1}{c|}{Definition} &  \multicolumn{1}{c|}{Value} \\ \hline
$\Omega$ & Domain & $(0,20)\times(0,8)$  \auth{($\si{mm}$)}\\ 
$ h $ & Diagonal cell diameter & test-dependent\\
$ G_C $ & Material toughness & $1.0\auth{\si{N}/\si{mm}}$ \\ 
$ E $ & Young's modulus & $1.0\auth{\si{kN}/\si{mm^2}}$ \\ 
$\mu$ & Lamé parameter &  $8.0\auth{\si{kN}/\si{mm^2}}$\\ 
$\lambda$ & Lamé parameter & $12.0\auth{\si{kN}/\si{mm^2}}$ \\ 
$ \nu $ & Poisson's ratio & $0.3$ \\ 
$ \varepsilon $ & Bandwidth of the initial crack & $2h$ \auth{($\si{mm}$)} \\
$ \kappa $ & Regularization parameter & \auth{$10^{-10}h$ ($\si{mm}$)}\\ 
$k_n$ & Time step size& \auth{$10^{-4}\si{s}$}\\
& \RA{Number of global refinements} & \RA{$2$}\\
& \RA{Number of local refinements} & \RA{$1$, $2$, $3$, $4$}\\
\RA{$\operatorname{TOL}_N$}& \RA{Tolerance outer Newton solver} & \RA{$10^{-7}$}  \\
& \RA{Tolerance inner linear solver} & \RA{$\lVert \Tilde{R}(U_h^{n,k})\rVert_2 10^{-8}$}  \\
\RA{$\operatorname{TOL}_{\operatorname{ItL}}$}& \RA{Tolerance ItL} & \RA{$10^{-3}$}  \\\hline
\end{tabular}
\caption{The setting of the material and numerical parameters used for the asymmetric three-point bending test.}
\label{table:params_astp}
\end{table}
\RA{
Two different final configurations are depicted in Figure \ref{fig:astp_vis_final_conf_less_ref} and Figure \ref{fig:astp_vis_final_conf_more_ref}. As it can be seen there, the crack path from the point, where the fracture exceeds the second hole 
depends on the mesh.
In \cite{MesBouKhon15} the authors also observe different crack paths depending on different solution approaches. However, as in \cite{AmGeraLoren15}, we are mainly interested in the simulation results before the fracture grows into the second hole. Thus, we always stop the simulation, once the crack reaches the second hole.
}

\begin{figure}[htbp!]
\begin{minipage}[b]{0.49\textwidth}
    \centering
    \includegraphics[scale = 0.21]{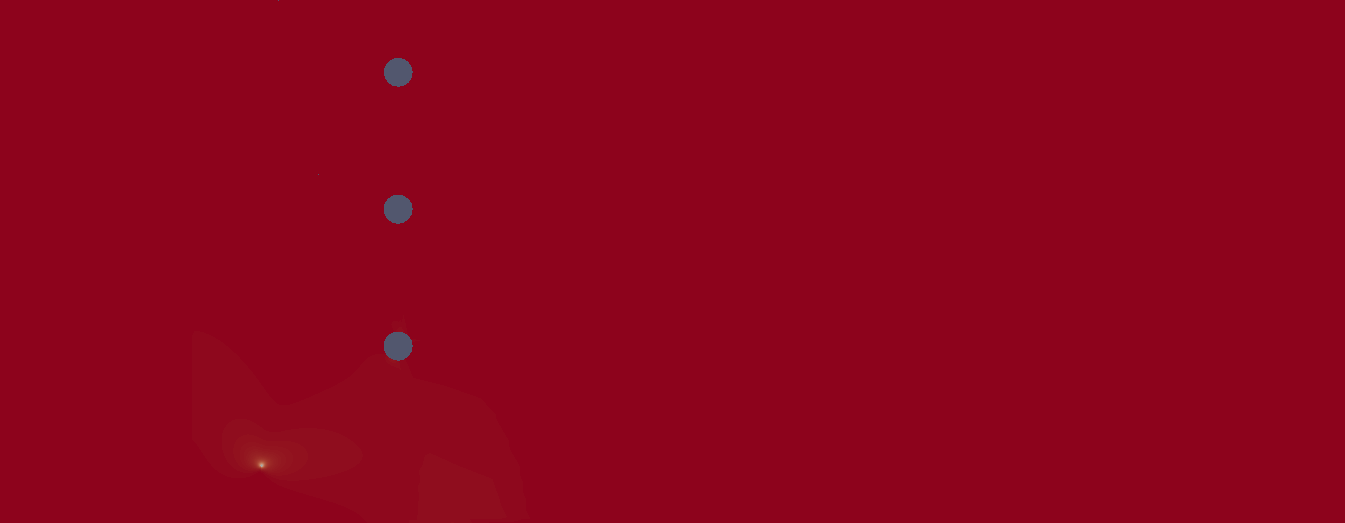}
    \caption{Visualization of the fracture asymmetric three-point bending test with $h=0.019$ after $146$ incremental steps with ItE. \RA{Red represents the fully intact part of the domain, blue the fully broken part and white stands for the transition zone.}}
    \label{fig:astp_vis_extra_146}
\end{minipage}
\hfill
\begin{minipage}[b]{0.49\textwidth}
    \centering
    \includegraphics[scale=0.21]{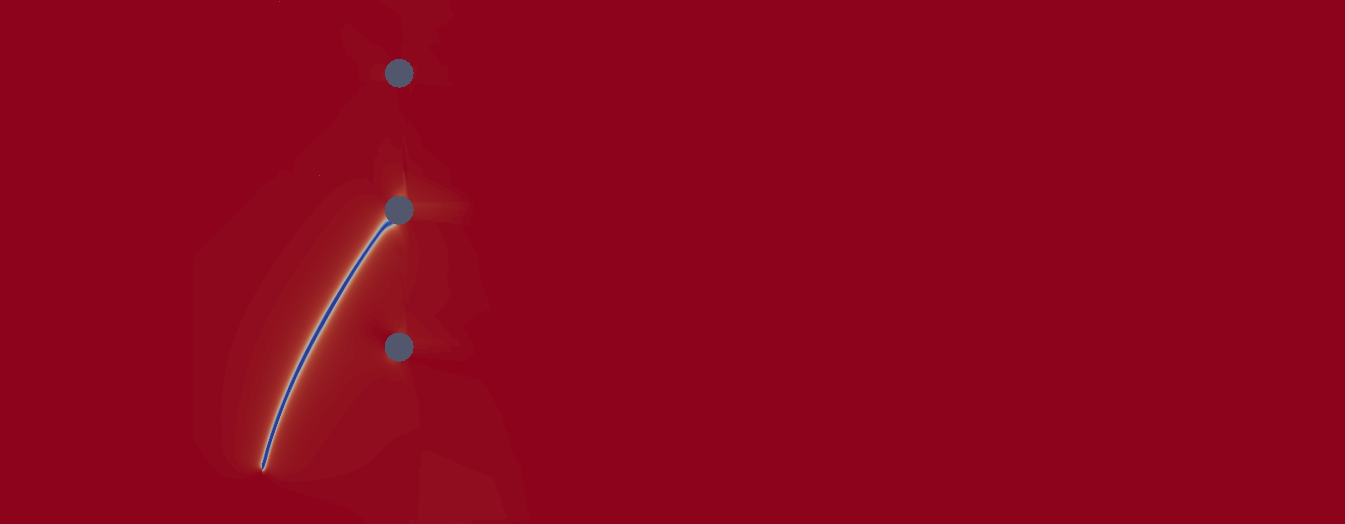}
    \caption{Visualization of the fracture asymmetric three-point bending test with $h=0.019$ after $147$ incremental steps with ItE. \RA{Red represents the fully intact part of the domain, blue the fully broken part and white stands for the transition zone.}}
    \label{fig:astp_vis_extra_147}
\end{minipage}
\end{figure}

In Figure \ref{plot_as_iter_0.156_astp} - Figure \ref{plot_load_disp_0.019_astp}, we compare the number of active set iterations and the crack energy for two different situations: no ItL and ItE \RA{with a tolerance of 
$\operatorname{TOL}_{\operatorname{ItL}} = 10^{-3}$ (see Algorithm \ref{alg:full_algorithm}) as $L^2$-difference. At the highest, around $160$ iterations on the linearizations where performed within one incremental step. Usually, within incremental steps without crack development, $2$ iterations are sufficient.} The crack energy can be computed via 
\begin{equation*}
    E_C = \frac{G_c}{2} \integral{\Omega}{}{\frac{(\varphi-1)^2}{\varepsilon}}{\Omega}.
\end{equation*}
We observe that with ItE, the material tears within one incremental step completely (from the tip of the initial fracture into the second hole, see Figure \ref{fig:astp_vis_extra_146} and Figure \ref{fig:astp_vis_extra_147}), whereas the crack evolves much slower without ItE. The former is the expected behaviour since the model is designed to represent brittle fractures. Additionally, without ItE, the time of full rupture differs much considering different refinement levels. This indicates that the influence of the mesh on the fracture development is less significant with ItE.

As before we want to give a short comment on the parallel performance for the mesh size 
$h=0.039$. Moreover, we investigate both without ItE and ItE. 
Without any iteration on the linearization and on $1$ core, 
the \RA{CPU} time for one incremental step is \RA{$678\mathrm{s}$} (approx. \RA{$11.3$} minutes), 
and on \RA{$16$} cores we observed a \RA{CPU} time of $59\mathrm{s}$.
With ItE, we have a \RA{CPU} time of 
\RA{$2207\mathrm{s}$} (approx. \RA{$36$ minutes})  on $1$ core 
and \RA{$174\mathrm{s}$} (approx. \RA{$2.9$} minutes) 
on \RA{$16$} cores \RA{for one incremental step (with $2$ ItE iterations)}.
This demonstrates clearly the balance of accuracy and efficiency whether ItE is used or not.

\begin{figure}[htbp!]
\begin{minipage}[b]{0.49\textwidth}
\scriptsize
\centering
\begin{tikzpicture}[scale = 1.0]
\begin{axis}[
    xlabel = time \auth{$[\si{s}]$}, 
    ymode =log,
    ylabel = $\#$ Active set iter\auth{ations},
    legend style = {at={(0.982,1.15)}},
    mark size=0.3pt,
    grid =major,
    y post scale = 0.9,
    ]
        \addplot[line width = 6pt,color=blue]
        table[col sep=space] {astp_as_iter_case_2_0.156.txt};
        \addlegendentry{$h=0.156$ $|$ no iteration on extrapolation}
        \addplot[line width = 6pt,color=orange, dashed]
        table[col sep=space] {astp_as_iter_case_2_0.156_extra.txt};
        \addlegendentry{$h=0.156$ $|$ with iteration on extrapolation}
\end{axis}
\end{tikzpicture}
\caption{Number of active set iterations ($y$-axis) with and without iteration on the extrapolation and $h=0.156$ depending on the time ($x$-axis) for the asymmetric three-point bending test.}
\label{plot_as_iter_0.156_astp}

\end{minipage}
\hfill
\begin{minipage}[b]{0.49\textwidth}
\scriptsize
\centering
\begin{tikzpicture}[scale = 1.0]
\begin{axis}[
    xlabel = time \auth{$[\si{s}]$}, 
    ylabel = Crack energy $E_C$ \auth{$[\si{\joule}]$},
    legend style = {at={(0.982,1.15)}},
    mark size=0.3pt,
    grid =major,
    y post scale = 0.9,
    ]
        \addplot[line width = 6pt,color=blue]
        table[col sep=space] {astp_crack_energy_case_2_no_extra_0.156.txt};
        \addlegendentry{$h=0.156$ $|$ no iteration on extrapolation}
        \addplot[line width = 6pt,color=orange, dashed]
        table[col sep=space] {astp_crack_energy_case_2_extra_0.156.txt};
        \addlegendentry{$h=0.156$ $|$ with iteration on extrapolation}
\end{axis}
\end{tikzpicture}
\caption{Visualization of the crack energy (y-axis) with and without iteration on the extrapolation and $h = 0.156$ depending on the time ($x$-axis) for the asymmetric three-point bending test.}
\label{plot_load_disp_0.156_astp}
\end{minipage}
\end{figure}

\begin{figure}[htbp!]
\begin{minipage}[b]{0.49\textwidth}
\scriptsize
\centering
\begin{tikzpicture}[scale = 1.0]
\begin{axis}[
    xlabel = time \auth{$[\si{s}]$}, 
    ymode = log, 
    ylabel = $\#$ Active set iter\auth{ations},
    legend style = {at={(0.982,1.15)}},
    mark size=0.3pt,
    grid =major,
    y post scale = 0.9,
    ]
        \addplot[line width = 6pt,color=blue]
        table[col sep=space] {astp_as_iter_case_2_0.078.txt};
        \addlegendentry{$h=0.078$ $|$ no iteration on extrapolation}
        \addplot[line width = 6pt,color=orange, dashed]
        table[col sep=space] {astp_as_iter_case_2_0.078_extra.txt};
        \addlegendentry{$h=0.078$ $|$ with iteration on extrapolation}
\end{axis}
\end{tikzpicture}
\caption{Number of active set iterations ($y$-axis) with and without iteration on the extrapolation and $h=0.078$ depending on the time ($x$-axis) for the asymmetric three-point bending test.}
\label{plot_as_iter_0.078_astp}

\end{minipage}
\hfill
\begin{minipage}[b]{0.49\textwidth}
\scriptsize
\centering
\begin{tikzpicture}[scale = 1.0]
\begin{axis}[
    xlabel = time \auth{$[\si{s}]$}, 
    ylabel = Crack energy $E_C$ \auth{$[\si{\joule}]$},
    legend style = {at={(0.982,1.15)}},
    mark size=0.3pt,
    grid =major,
    y post scale = 0.9,
    ]
        \addplot[line width = 6pt,color=blue]
        table[col sep=space] {astp_crack_energy_case_2_no_extra_0.078.txt};
        \addlegendentry{$h=0.078$ $|$ no iteration on extrapolation}
        \addplot[line width = 6pt,color=orange, dashed]
        table[col sep=space] {astp_crack_energy_case_2_extra_0.078.txt};
        \addlegendentry{$h=0.078$ $|$ with iteration on extrapolation}
\end{axis}
\end{tikzpicture}
\caption{Visualization of the crack energy (y-axis) with and without iteration on the extrapolation and $h = 0.078$ depending on the time ($x$-axis) for the asymmetric three-point bending test.}
\label{plot_load_disp_0.078_astp}
\end{minipage}
\end{figure}

\begin{figure}[htbp!]
\begin{minipage}[b]{0.49\textwidth}
\scriptsize
\centering
\begin{tikzpicture}[scale = 1.0]
\begin{axis}[
    xlabel = time \auth{$[\si{s}]$}, 
    ymode = log, 
    ylabel = $\#$ Active set iter\auth{ations},
    legend style = {at={(0.982,1.15)}},
    mark size=0.3pt,
    grid =major,
    y post scale = 0.9,
    ]
        \addplot[line width = 6pt,color=blue]
        table[col sep=space] {astp_as_iter_case_2_0.039.txt};
        \addlegendentry{$h=0.039$ $|$ no iteration on extrapolation}
        \addplot[line width = 6pt,color=orange, dashed]
        table[col sep=space] {astp_as_iter_case_2_0.039_extra.txt};
        \addlegendentry{$h=0.039$ $|$ with iteration on extrapolation}
        
\end{axis}
\end{tikzpicture}
\caption{Number of active set iterations ($y$-axis) with and without iteration on the extrapolation and $h=0.039$ depending on the time ($x$-axis) for the asymmetric three-point bending test.}
\label{plot_as_iter_0.039_astp}

\end{minipage}
\hfill
\begin{minipage}[b]{0.49\textwidth}
\scriptsize
\centering
\begin{tikzpicture}[scale = 1.0]
\begin{axis}[
    xlabel = time \auth{$[\si{s}]$}, 
    ylabel = Crack energy $E_C$ \auth{$[\si{\joule}]$},
    legend style = {at={(0.982,1.15)}},
    mark size=0.3pt,
    grid =major,
    y post scale = 0.9,
    ]
        \addplot[line width = 6pt,color=blue]
        table[col sep=space] {astp_crack_energy_case_2_no_extra_0.039.txt};
        \addlegendentry{$h=0.039$ $|$ no iteration on extrapolation}
        \addplot[line width = 6pt,color=orange, dashed]
        table[col sep=space] {astp_crack_energy_case_2_extra_0.039.txt};
        \addlegendentry{$h=0.039$ $|$ with iteration on extrapolation}
\end{axis}
\end{tikzpicture}
\caption{Visualization of the crack energy (y-axis) with and without iteration on the extrapolation and $h = 0.039$ depending on the time ($x$-axis) for the asymmetric three-point bending test.}
\label{plot_load_disp_0.039_astp}
\end{minipage}
\end{figure}

\begin{figure}[htbp!]
\begin{minipage}[b]{0.49\textwidth}
\scriptsize
\centering
\begin{tikzpicture}[scale = 1.0]
\begin{axis}[
    xlabel = time \auth{$[\si{s}]$}, 
    ymode = log,
    ylabel = $\#$ Active set iter\auth{ations},
    legend style = {at={(0.982,1.15)}},
    mark size=0.3pt,
    grid =major,
    y post scale = 0.9,
    ]
        \addplot[line width = 6pt,color=blue, dashed]
        table[col sep=space] {astp_as_iter_case_2_0.019.txt};
        \addlegendentry{$h=0.019$ $|$ no iteration on extrapolation}
        \addplot[line width = 6pt,color=orange, dotted]
        table[col sep=space] {astp_as_iter_case_2_0.019_extra.txt};
        \addlegendentry{$h=0.019$ $|$ with iteration on extrapolation}
        
\end{axis}
\end{tikzpicture}
\caption{Number of active set iterations ($y$-axis) with and without iteration on the extrapolation and $h=0.019$ depending on the time ($x$-axis) for the asymmetric three-point bending test.}
\label{plot_as_iter_0.019_astp}

\end{minipage}
\hfill
\begin{minipage}[b]{0.49\textwidth}
\scriptsize
\centering
\begin{tikzpicture}[scale = 1.0]
\begin{axis}[
    xlabel = time \auth{$[\si{s}]$}, 
    ylabel = Crack energy $E_C$ \auth{$[\si{\joule}]$},
    legend style = {at={(0.982,1.15)}},
    mark size=0.3pt,
    grid =major,
    y post scale = 0.9,
    ]
        \addplot[line width = 6pt,color=blue]
        table[col sep=space] {astp_crack_energy_case_2_no_extra_0.019.txt};
        \addlegendentry{$h=0.019$ $|$ no iteration on extrapolation}
        \addplot[line width = 6pt,color=orange, dashed]
        table[col sep=space] {astp_crack_energy_case_2_extra_0.019.txt};
        \addlegendentry{$h=0.019$ $|$ with iteration on extrapolation}
 \end{axis}
\end{tikzpicture}
\caption{Visualization of the crack energy (y-axis) with and without iteration on the extrapolation and $h = 0.019$ depending on the time ($x$-axis) for the asymmetric three-point bending test.}
\label{plot_load_disp_0.019_astp}
\end{minipage}
\end{figure}

\newpage
\subsection{\RA{L-shaped panel test}}
\begin{figure}[htbp!]
\centering
\begin{tikzpicture}[xscale=0.01,yscale=0.01]
\fill[gray!30] (0,0) -- (250,) -- (250,250) -- (500,250) -- (500,500) -- (0,500) -- (0,0);
\fill [gray!50] (500,500) rectangle (400,250);
\draw[] (0,0) -- (250,) -- (250,250) -- (500,250) -- (500,500) -- (0,500) -- (0,0);
\draw[<->] (0,-40) -- (250,-40);
\draw (15,0) -- (0,-20);
\draw (35,0) -- (20,-20);
\draw (55,0) -- (40,-20);
\draw (75,0) -- (60,-20);
\draw (95,0) -- (80,-20);
\draw (115,0) -- (100,-20);
\draw (135,0) -- (120,-20);
\draw (155,0) -- (140,-20);
\draw (175,0) -- (160,-20);
\draw (195,0) -- (180,-20);
\draw (215,0) -- (200,-20);
\draw (235,0) -- (220,-20);
\node[anchor = north] at (125,-50) {$250\,\mathrm{mm}$};
\node[anchor = north] at (125,-95) {$\partial \Omega_{\operatorname{measured}}$};
\draw[<->] (-40,0) -- (-40,500);
\node[anchor = east] at (-50,250) {$500\,\mathrm{mm}$};
\draw[<->] (0,540) -- (500,540);
\node[anchor = south] at (250,550) {$500\,\mathrm{mm}$};
\draw[<->] (540,500) -- (540,250);
\node[anchor = west] at (550,375) {$250\,\mathrm{mm}$};
\draw[blue] (500,250) -- (470,250);
\draw[blue, ->] (475,210) -- (475,240);
\draw[blue, <-] (495,210) -- (495,240);
\node[anchor = north] at (485,210) {$u_y$};
\node[anchor = north] at (485,180) {$30\,\mathrm{mm}$};
\draw[<->] (400,375) -- (500,375);
\node[anchor = south] at (430,385) {$100\,\mathrm{mm}$};
\end{tikzpicture}
 \caption{\RA{Visualization of the geometry of the L-shaped panel test.}}\label{fig:l_shaped_geo}
\end{figure}
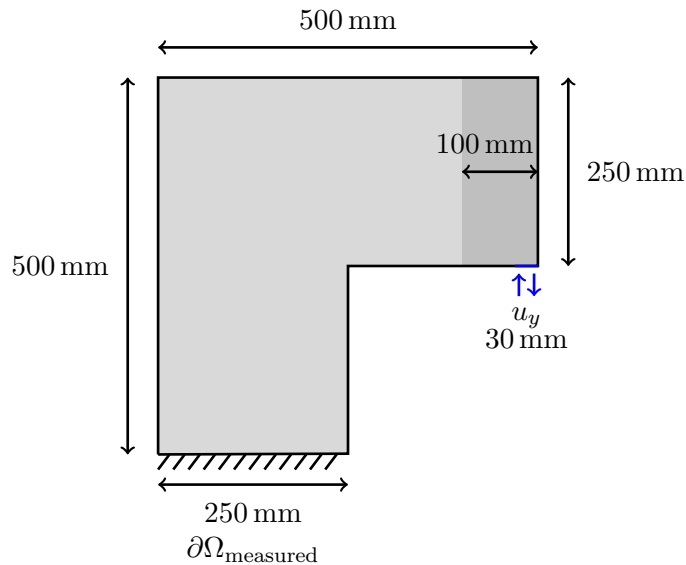

\RA{The L-shaped panel test is a well-known test 
from mechanics, 
which was originally developed by Winkler \cite{winkler2001traglastuntersuchungen} 
to analyze possible crack behaviour of concrete under force
experimentally as well as numerically. Further numerical simulations in connection to
variational/phase-field solution approaches are performed in \cite{AmGeraLoren15, GeLo16, MesBouKhon15, Wi17_SISC,MaWaWiWo20, mang2022phasefield}. In Figure \ref{fig:l_shaped_geo}, the geometry $\Omega$ of the L-shaped panel test is depicted.
There exist two configurations of the loading boundary conditions: first monotone loading, and second, a cyclic loading test. In the 
following, we have the latter in mind.
}

\RA{
Specifically, we apply a cyclic displacement as non-homogeneous time-dependent Dirichlet conditions \cite{AmGeraLoren15} at a strip of $30\si{mm}$ length on the right corner of the domain:
\begin{alignat}{2}
    u_y &= t, \quad &&0.0\si{s} \leq t \leq 0.3\si{s} \\
    u_y &= 0.6-t, \quad &&0.3\si{s} \leq t \leq 0.8\si{s} \\
    u_y &= t-1.0, \quad &&0.8\si{s} \leq t \leq 2.0\si{s},
\end{alignat}
where $t$ represents the simulation time. On the lower boundary $\partial \Omega_{\operatorname{measure}}$, 
we fix the displacement in 
$x$- and $y$-direction with homogeneous Dirichlet conditions as $u_x = u_y=0$. 
On the other boundaries, we apply homogeneous Neumann conditions.
As done in \cite{MesBouKhon15}, in order to avoid the development of 
nonphysical fractures around the right edge of the specimen, 
the phase field is constrained with $\varphi = 1$ for $x > 400\si{mm}$.
As quantities of interest, we observe the number of Newton iterations and the load on the bottom boundary $\partial \Omega_{\operatorname{measure}}$, computed via
\begin{align}
 (F_x,F_y)\coloneqq \int_{\partial \Omega_{\operatorname{measure}}} \sigma(u_h)\cdot \eta\,\mathrm{d}s,\label{loading_l_shaped}
\end{align}
where $\eta$ is the unit normal vector.
The parameters for this test are displayed in Table \ref{table:params_l_shaped}}. 
\begin{table}[htbp!]
\renewcommand{\arraystretch}{1.2}
\scriptsize
\centering
\RA{
\begin{tabular}{|l|l|c|}\hline
\multicolumn{1}{|c|}{Parameter} & \multicolumn{1}{c|}{Definition} &  \multicolumn{1}{c|}{Value} \\ \hline
$\Omega$ & Domain & $\left((0,0)\times(250,500)\right)\cup \left((250,250)\times (500,500)\right)$ ($\si{mm}$)  \\ 
$ h $ & Diagonal cell diameter & test-dependent\\
$ G_C $ & Material toughness & $8.9\cdot 10^{-2}\si{N}/\si{mm}$ \\ 
$ E $ & Young's modulus & $10.677333\si{kN}/\si{mm^2}$ \\ 
$\mu$ & Lamé parameter &  $10.95\si{kN}/\si{mm^2}$\\ 
$\lambda$ & Lamé parameter & $6.16\si{kN}/\si{mm^2}$ \\ 
$ \nu $ & Poisson's ratio & $0.3$ \\ 
$ \varepsilon $ & Bandwidth of the initial crack & $2h$ ($\si{mm}$) \\
$ \kappa $ & Regularization parameter & $10^{-10}h$ ($\si{mm}$)\\ 
$k_n$ & Time step size& $10^{-3}\si{s}$\\
& Number of global refinements & $2$, $3$, $4$, $5$\\
& Number of local refinements & $0$\\
$\operatorname{TOL}_N$& Tolerance outer Newton solver & $10^{-7}$  \\
& Tolerance inner linear solver & $\lVert \Tilde{R}(U_h^{n,k})\rVert_2 10^{-8}$  \\
\RA{$\operatorname{TOL}_{\operatorname{ItL}}$}& \RA{Tolerance ItL} & \RA{$10^{-1}$}  \\\hline
\end{tabular}
}
\caption{\RA{The setting of the material and numerical parameters used for the L-shaped panel test.}}
\label{table:params_l_shaped}
\end{table}

\begin{figure}
    \centering
    \includegraphics[scale = 0.20]{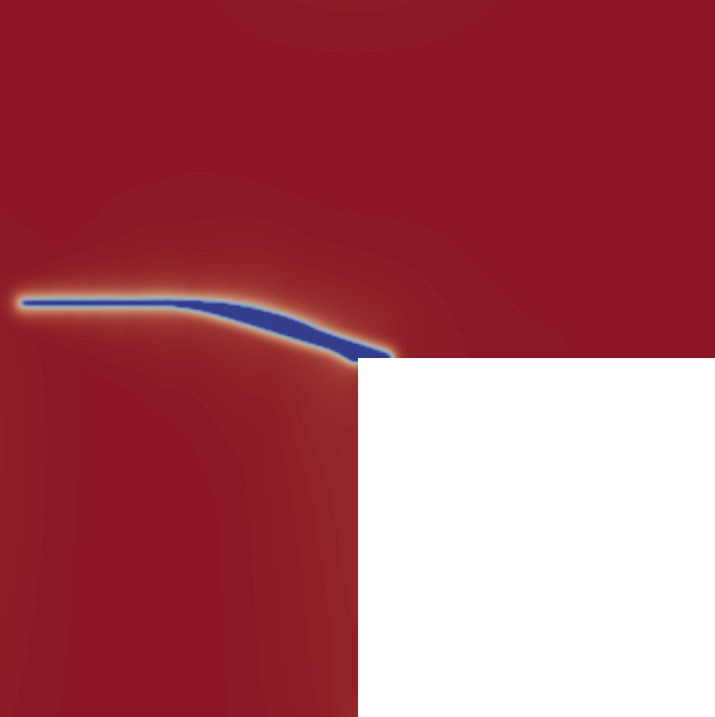}
    \caption{\RA{Visualization of the final configuration of the L-shaped panel test with $h=1.822$.} \RA{Red represents the fully intact part of the domain, blue the fully broken part and white stands for the transition zone.}}
    \label{fig:l_shaped_final_conf}
\end{figure}

\RA{A typical final configuration on a globally refined mesh with diagonal element diameter $h = 1.822$ is displayed in Figure \ref{fig:l_shaped_final_conf}. Figure \ref{plot_as_iter_l_shaped_h_14.577}-\ref{plot_load_disp_l_shaped_h_1.822} visualize the number of active set iterations (as before plus Newton, i.e., combined Newton iterations are given) per incremental step using \textbf{Case } 2 from Section \ref{subsubsec:4_cases} and the load displacement curves with no ItL as well as ItE with a tolerance of 
$\operatorname{TOL}_{\operatorname{ItL}} = 10^{-1}$ (see Algorithm \ref{alg:full_algorithm}) as $L^2$-difference and for four different global refinement levels. On each refinement level, we observe the benefit of ItE as the load decreases faster than with no ItL, i.e. the fracture evolves faster. However, the additional computational cost is also obvious as in the incremental steps of crack evolvement the number of active set iterations is much higher than with no ItL.}

\RA{For the parallel performance we compare, as in the previous examples, the CPU time for one incremental step on a $5$ times globally refined mesh with $h = 1.82217$ ($232323$ degrees of freedom). With ItE and on $1$ core, the CPU time for one incremental step (and $2$ iterations on the extrapolation) is $2094\si{s}$ (approx. $34.9$ minutes), whereas on $16$ cores it is $154\si{s}$ (approx. $2.6$ minutes). With no ItL and on $1$ core, $1109\si{s}$ (approx. $18.5$ minutes) are needed and $76\si{s}$ (approx. $1.3$ minutes) on $16$ cores.}
\RA{
\begin{figure}[htbp!]
\begin{minipage}[b]{0.49\textwidth}
\scriptsize
\centering
\begin{tikzpicture}[scale = 1.0]
\begin{axis}[
    xlabel = time $\lbrack\si{s}\rbrack$, 
    ymode = log,
    ylabel = $\#$ Active set iterations,
    legend style = {at={(0.75,1.3)}},
    mark size=0.3pt,
    grid =major,
    y post scale = 0.9,
    ]
        \addplot[line width = 6pt,color=blue]
        table[col sep=space] {as_iter_no_ItL_2_0.txt};
        \addlegendentry{No ItL $|$ $h=14.577$}
        \addplot[line width = 6pt, orange, dashed]
        table[col sep=space] {as_iter_ItE_2_0.txt};
        \addlegendentry{ItE $|$ $h=14.577$}
\end{axis}
\end{tikzpicture}
\caption{\RA{Visualization of the number of active set iterations ($y$-axis) with and without iteration on the extrapolation (ItE) and $h=14.577$ depending on the time ($x$-axis) for the L-shaped panel test.}}
\label{plot_as_iter_l_shaped_h_14.577}
\end{minipage}
\hfill
\begin{minipage}[b]{0.49\textwidth}
\scriptsize
\centering
\begin{tikzpicture}[scale = 1.0]
\begin{axis}[
    xlabel = displacement $u$ $\lbrack\si{mm}\rbrack$, 
    ylabel = Load $F_y$ $\lbrack\si{\newton}\rbrack$,
    legend style = {at={(0.75,1.3)}},
    mark size=0.3pt,
    grid =major,
    y post scale = 0.9,
    ]
        \addplot[line width = 6pt,color=blue]
        table[col sep=space] {load_y_no_ItL_2_0.txt};
        \addlegendentry{No ItL $|$ $h=14.577$}
        \addplot[line width = 6pt, orange, dashed]
        table[col sep=space] {load_y_ItE_2_0.txt};
        \addlegendentry{ItE $|$ $h=14.577$}
\end{axis}
\end{tikzpicture}
\caption{\RA{Visualization of the $y$-load (y-axis) with and without iteration on the extrapolation (ItE) and $h = 14.577$ depending on the displacement ($x$-axis) for the L-shaped panel test.}}
\label{plot_load_disp_l_shaped_h_14.577}
\end{minipage}
\end{figure}
}

\RA{
\begin{figure}[htbp!]
\begin{minipage}[b]{0.49\textwidth}
\scriptsize
\centering
\begin{tikzpicture}[scale = 1.0]
\begin{axis}[
    xlabel = time $\lbrack\si{s}\rbrack$, 
    ymode = log,
    ylabel = $\#$ Active set iterations,
    legend style = {at={(0.75,1.3)}},
    mark size=0.3pt,
    grid =major,
    y post scale = 0.9,
    ]
        \addplot[line width = 6pt,color=blue]
        table[col sep=space] {as_iter_no_ItL_3_0.txt};
        \addlegendentry{No ItL $|$ $h=7.288$}
        \addplot[line width = 6pt, orange, dashed]
        table[col sep=space] {as_iter_ItE_3_0.txt};
        \addlegendentry{ItE $|$ $h=7.288$}
\end{axis}
\end{tikzpicture}
\caption{\RA{Visualization of the number of active set iterations ($y$-axis) with and without iteration on the extrapolation (ItE) and $h=7.288$ depending on the time ($x$-axis) for the L-shaped panel test.}}
\label{plot_as_iter_l_shaped_h_7.288}
\end{minipage}
\hfill
\begin{minipage}[b]{0.49\textwidth}
\scriptsize
\centering
\begin{tikzpicture}[scale = 1.0]
\begin{axis}[
    xlabel = displacement $u$ $\lbrack\si{mm}\rbrack$, 
    ylabel = Load $F_y$ $\lbrack\si{\newton}\rbrack$,
    legend style = {at={(0.75,1.3)}},
    mark size=0.3pt,
    grid =major,
    y post scale = 0.9,
    ]
        \addplot[line width = 6pt,color=blue]
        table[col sep=space] {load_y_no_ItL_3_0.txt};
        \addlegendentry{No ItL $|$ $h=7.288$}
        \addplot[line width = 6pt, orange, dashed]
        table[col sep=space] {load_y_ItE_3_0.txt};
        \addlegendentry{ItE $|$ $h=7.288$}
\end{axis}
\end{tikzpicture}
\caption{\RA{Visualization of the $y$-load (y-axis) with and without iteration on the extrapolation (ItE) and $h = 7.288$ depending on the displacement ($x$-axis) for the L-shaped panel test.}}
\label{plot_load_disp_l_shaped_h_7.288}
\end{minipage}
\end{figure}
}

\RA{
\begin{figure}[htbp!]
\begin{minipage}[b]{0.49\textwidth}
\scriptsize
\centering
\begin{tikzpicture}[scale = 1.0]
\begin{axis}[
    xlabel = time $\lbrack\si{s}\rbrack$, 
    ymode = log,
    ylabel = $\#$ Active set iterations,
    legend style = {at={(0.75,1.3)}},
    mark size=0.3pt,
    grid =major,
    y post scale = 0.9,
    ]
        \addplot[line width = 6pt,color=blue]
        table[col sep=space] {as_iter_no_ItL_4_0.txt};
        \addlegendentry{No ItL $|$ $h=3.644$}
        \addplot[line width = 6pt, orange, dashed]
        table[col sep=space] {as_iter_ItE_4_0.txt};
        \addlegendentry{ItE $|$ $h=3.644$}
\end{axis}
\end{tikzpicture}
\caption{\RA{Visualization of the number of active set iterations ($y$-axis) with and without iteration on the extrapolation (ItE) and $h=3.644$ depending on the time ($x$-axis) for the L-shaped panel test.}}
\label{plot_as_iter_l_shaped_h_3.644}
\end{minipage}
\hfill
\begin{minipage}[b]{0.49\textwidth}
\scriptsize
\centering
\begin{tikzpicture}[scale = 1.0]
\begin{axis}[
    xlabel = displacement $u$ $\lbrack\si{mm}\rbrack$, 
    ylabel = Load $F_y$ $\lbrack\si{\newton}\rbrack$,
    legend style = {at={(0.75,1.3)}},
    mark size=0.3pt,
    grid =major,
    y post scale = 0.9,
    ]
        \addplot[line width = 6pt,color=blue]
        table[col sep=space] {load_y_no_ItL_4_0.txt};
        \addlegendentry{No ItL $|$ $h=3.644$}
        \addplot[line width = 6pt, orange, dashed]
        table[col sep=space] {load_y_ItE_4_0.txt};
        \addlegendentry{ItE $|$ $h=3.644$}
\end{axis}
\end{tikzpicture}
\caption{\RA{Visualization of the $y$-load (y-axis) with and without iteration on the extrapolation (ItE) and $h = 3.644$ depending on the displacement ($x$-axis) for the L-shaped panel test.}}
\label{plot_load_disp_l_shaped_h_3.644}
\end{minipage}
\end{figure}
}

\RA{
\begin{figure}[htbp!]
\begin{minipage}[b]{0.49\textwidth}
\scriptsize
\centering
\begin{tikzpicture}[scale = 1.0]
\begin{axis}[
    xlabel = time $\lbrack\si{s}\rbrack$, 
    ymode = log,
    ylabel = $\#$ Active set iterations,
    legend style = {at={(0.75,1.3)}},
    mark size=0.3pt,
    grid =major,
    y post scale = 0.9,
    ]
        \addplot[line width = 6pt,color=blue]
        table[col sep=space] {as_iter_no_ItL_5_0.txt};
        \addlegendentry{No ItL $|$ $h=1.822$}
        \addplot[line width = 6pt, orange, dashed]
        table[col sep=space] {as_iter_ItE_5_0.txt};
        \addlegendentry{ItE $|$ $h=1.822$}
\end{axis}
\end{tikzpicture}
\caption{\RA{Visualization of the number of active set iterations ($y$-axis) with and without iteration on the extrapolation (ItE) and $h=1.822$ depending on the time ($x$-axis) for the L-shaped panel test.}}
\label{plot_as_iter_l_shaped_h_1.822}
\end{minipage}
\hfill
\begin{minipage}[b]{0.49\textwidth}
\scriptsize
\centering
\begin{tikzpicture}[scale = 1.0]
\begin{axis}[
    xlabel = displacement $u$ $\lbrack\si{mm}\rbrack$, 
    ylabel = Load $F_y$ $\lbrack\si{\newton}\rbrack$,
    legend style = {at={(0.75,1.3)}},
    mark size=0.3pt,
    grid =major,
    y post scale = 0.9,
    ]
        \addplot[line width = 6pt,color=blue]
        table[col sep=space] {load_y_no_ItL_5_0.txt};
        \addlegendentry{No ItL $|$ $h=1.822$}
        \addplot[line width = 6pt, orange, dashed]
        table[col sep=space] {load_y_ItE_5_0.txt};
        \addlegendentry{ItE $|$ $h=1.822$}
\end{axis}
\end{tikzpicture}
\caption{\RA{Visualization of the $y$-load (y-axis) with and without iteration on the extrapolation (ItE) and $h = 1.822$ depending on the displacement ($x$-axis) for the L-shaped panel test.}}
\label{plot_load_disp_l_shaped_h_1.822}
\end{minipage}
\end{figure}
}
\newpage
\subsection{Single edge notched shear test}\label{subsec:sens}
\begin{figure}[htbp!]
\centering
 \begin{tikzpicture}[xscale=0.85,yscale=0.85]
\draw[fill=gray!30] (0,0) -- (0,5) -- (5,5) -- (5,0) -- (0,0);
\draw[red] (2.5,2.5) -- (5,2.5);
\node[red] at (3.725,2.15) {slit};
\draw[<->,black] (2.5,2.75) -- (5,2.75);
\node[anchor = south] at (3.75,2.75) {\auth{$0.5\,\mathrm{mm}$}};
\draw (5,0) -- (5,5);
\draw[blue] (5,5) -- (0,5);
\draw[-] (0,5) -- (0,0);
\draw[->,blue] (0.7,5.5) -- (0.1,5.5);
\node at (0.4,5.25) {$u_x$};
\draw[<->] (-0.25,0) -- (-0.25,5);
\node[anchor = east] at (-0.25,2.5) {\auth{$1\,\mathrm{mm}$}};
\draw[<->] (0,-0.5) -- (5,-0.5);
\node at (5.15,-0.5) {$x$};
\node at (-0.5,5.15) {$y$};
\node[anchor = north] at (2.5,-0.5) {\auth{$1\,\mathrm{mm}$}};

\draw[<->] (5.25,0) -- (5.25,2.5);
\node[anchor = west] at (5.25,1.25) {\auth{$0.5\,\mathrm{mm}$}};

\draw (0.1,0) -- (-0.05,-0.2);
\draw (0.3,0) -- (0.15,-0.2);
\draw (0.5,0) -- (0.35,-0.2);
\draw (0.7,0) -- (0.55,-0.2);
\draw (0.9,0) -- (0.75,-0.2);
\draw (1.1,0) -- (0.95,-0.2);
\draw (1.3,0) -- (1.15,-0.2);
\draw (1.5,0) -- (1.35,-0.2);
\draw (1.7,0) -- (1.55,-0.2);
\draw (1.9,0) -- (1.75,-0.2);
\draw (2.1,0) -- (1.95,-0.2);
\draw (2.3,0) -- (2.15,-0.2);
\draw (2.5,0) -- (2.35,-0.2);
\draw (2.7,0) -- (2.55,-0.2);
\draw (2.9,0) -- (2.75,-0.2);
\draw (3.1,0) -- (2.95,-0.2);
\draw (3.3,0) -- (3.15,-0.2);
\draw (3.5,0) -- (3.35,-0.2);
\draw (3.7,0) -- (3.55,-0.2);
\draw (3.9,0) -- (3.75,-0.2);
\draw (4.1,0) -- (3.95,-0.2);
\draw (4.3,0) -- (4.15,-0.2);
\draw (4.5,0) -- (4.35,-0.2);
\draw (4.7,0) -- (4.55,-0.2);
\draw (4.9,0) -- (4.75,-0.2);
 \end{tikzpicture}
 \caption{Visualization of the geometry of the single edge notched shear test.}\label{fig:shear_geo}
\end{figure}
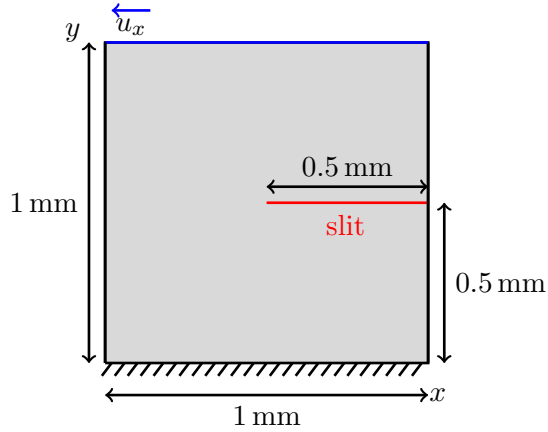

\begin{figure}[htbp!]
    \centering
    \includegraphics[scale=0.16]{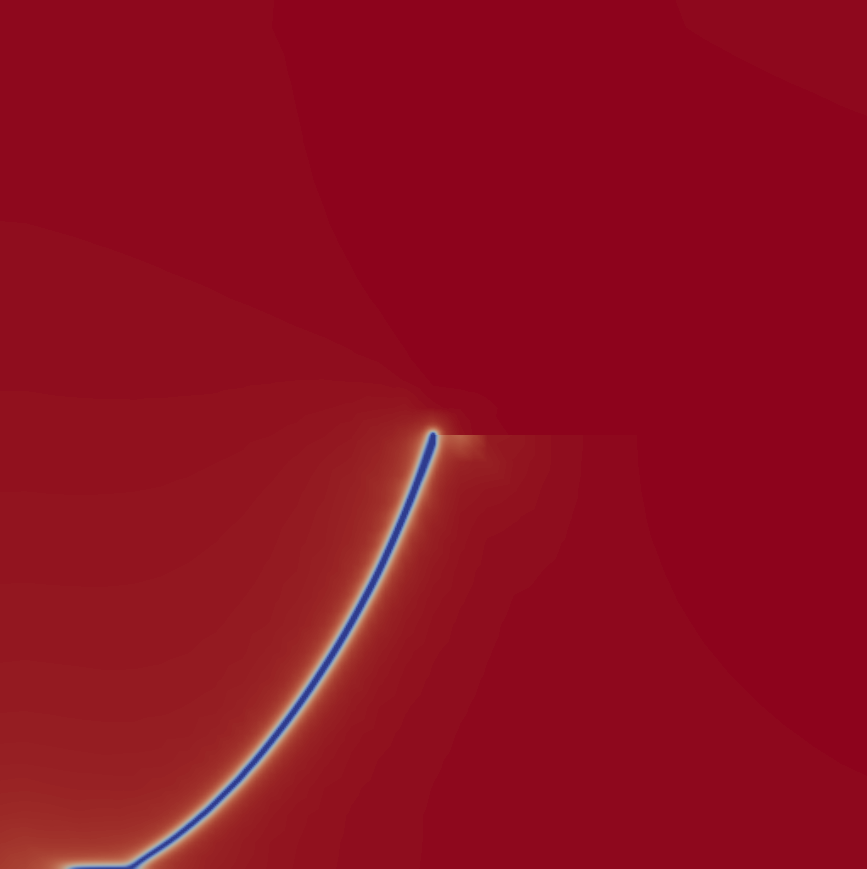}
    \caption{\RA{Visualization of the fracture for the SENS test with $h=0.0027$ in the final configuration after $135$ incremental steps with ItE.} \RA{Red represents the fully intact part of the domain, blue the fully broken part and white stands for the transition zone.}}
    \label{fig:sens_vis_final_conf}
\end{figure}

In the last example, we consider the single-edge notched shear (SENS), see for instance \cite{MieWelHof10a,miehe2010phase}. A visualization and explanation is displayed in Figure~\ref{fig:shear_geo}.
The domain $\Omega$ is a two dimensional square of \auth{$1\,\mathrm{mm}$} length with a given crack (called geometrical slit) on the right side at $5\,\mathrm{mm}$ tending to the midpoint of the square. On the bottom boundary, the square is fixed, and on the top boundary, a given displacement in the $x$-direction pulls to the left. The boundary conditions can be found for instance in \cite{Wi17_CMAME}: On the left and right sides, the boundaries are defined to be traction-free (homogeneous Neumann conditions). The bottom boundary is fixed via $u_x = u_y = 0\,\mathrm{mm}$. On the top boundary, it holds $u_y = 0\,\mathrm{mm}$ and in the $x$-direction we determine a time-dependent non-homogeneous Dirichlet condition: $u_x = t_n \cdot 1\,\mathrm{mm/s}$, where $t_n$ is the $n$-th incremental step. The end time $T$ is the incremental step once the specimen is broken. The parameters for this test are given in Table \ref{table_param_2}. \RA{A standard final configuration of this test example is depicted in Figure \ref{fig:sens_vis_final_conf}.} 

As a quantity of interest, we evaluate the load functions on the top boundary $\partial \Omega_{\operatorname{top}}$ computed via
\auth{
\begin{align}
 (F_x,F_y)\coloneqq \int_{\partial \Omega_{\operatorname{top}}} \sigma(u_h)\cdot \eta\,\mathrm{d}s,\label{loading1}
\end{align}
}
with the stress tensor $\sigma(u_h)$ depending on the discrete solution variable $u_h$ and the outer normal vector $\eta$. Within the single-edge notched shear test, we are interested in the evaluation of $F_x$.
\begin{table}[htbp!]
\renewcommand{\arraystretch}{1.2}
\scriptsize
\centering
\begin{tabular}{|l|l|c|}\hline
\multicolumn{1}{|c|}{Parameter} & \multicolumn{1}{c|}{Definition} &  \multicolumn{1}{c|}{Value} \\ \hline
$\Omega$ & Domain & $(0,1)^2$ \auth{($\si{mm}$)} \\ 
$ h $ & Diagonal cell diameter & test-dependent\\
$ G_C $ & Material toughness & $2.7\auth{\si{N}/\si{mm}}$ \\ 
$\mu$ & Lamé parameter &  $80.77\auth{\si{kN}/\si{mm^2}}$\\ 
$\lambda$ & Lamé parameter & $121.15\auth{\si{kN}/\si{mm^2}}$ \\ 
$ \nu $ & Poisson's ratio & $0.2$ \\ 
$ \varepsilon $ & Bandwidth of the initial crack & $2h$ \auth{($\si{mm}$)} \\
$ \kappa $ & Regularization parameter & \auth{$10^{-10}h$ ($\si{mm}$)}\\ 
$k_n$ & Time step size& $\auth{10^{-4}\si{s}}$\\
& \RA{Number of global refinements} & \RA{$2$}\\
& \RA{Number of local refinements} & \RA{$1$, $2$, $3$, $4$}\\
\RA{$\operatorname{TOL}_N$}& \RA{Tolerance outer Newton solver} & \RA{$10^{-7}$}  \\
& \RA{Tolerance inner linear solver} & \RA{$\lVert \Tilde{R}(U_h^{n,k})\rVert_2 10^{-8}$}  \\
\RA{$\operatorname{TOL}_{\operatorname{ItL}}$}& \RA{Tolerance ItL} & \RA{$10^{-1}$}  \\\hline

\end{tabular}
\caption{The setting of the material and numerical parameters used for the SENS-test.}
\label{table_param_2}
\end{table}

Similar to the asymmetric three-point bending test, we want to investigate the advantages and possible drawbacks of an iteration into the monolithic limit using ItL. As before, we compare the crack energy and the corresponding Newton active set iterations for two different approaches: In one situation, we only compute the linearization and solve the system with the Newton method once per incremental step. In the second situation, we iterate, as described before, until the $L^2$-norm of the difference between two consecutive solutions is small enough up to a given tolerance. 
Since we observe non-physical behaviour 
(Figure \ref{plot_load_disp_extra_iter_0.011_ms})
when using ItE\footnote{We are aware of \cite{AmGeraLoren15} who found that the Miehe et al.
splitting \cite{MieWelHof10a} is not a good choice here. But this does not alter our solver 
investigations we have primarily in mind.}, we also consider ItOTS as a third approach. We compare the total number of Newton active set iterations per incremental step and the $x$-load. For both ItE and ItOTS, we choose a tolerance of 
\RA{$\operatorname{TOL}_{\operatorname{ItL}} = 10^{-1}$ (see Algorithm \ref{alg:full_algorithm})} as $L^2$-difference. The results are visualized in Figure \ref{plot_load_disp_extra_iter_0.022_ms} - Figure \ref{plot_load_disp_extra_iter_0.0027_ms}. \RA{We also examine varying ItL tolerances and Newton tolerances, but found that both do not influence the solution. For ItL, a smaller tolerance only leads to a larger computation time since more iterations are necessary, but the solution does not differ significantly. For the Newton tolerance, the residual norm converges faster than the active set, such that the final residual norm is of magnitude around $10^{-10}$ no matter which tolerance is handed over.} 

We observe that on all refinement levels that fracture evolves much faster with ItE than with ItOTS or no ItL. This confirms findings 
made in \cite{Wi17_SISC}[Fig. 3] in which the extrapolated scheme (no ItL) is compared to a fully monolithic scheme.
However with ItE, for $h=0.022,\, 0.011,\, 0.0055$, the load shows unphysical behaviour: at a certain point, the fracture stops growing and the load increases again until it reaches a local peak. Then, the crack continues to propagate until the material is fully ruptured. From a physical perspective, we expect the fracture to fully evolve within a few incremental steps once it starts to evolve. This behaviour is less significant, when using ItOTS or no ItL. For the number of Newton active set iterations, no ItL needs less iterations than ItE and ItOTS. This is clear since without ItL, in each incremental step the system is only solved once per predictor-corrector refinement step, whereas with ItE and ItOTS, the system may be solved many times per predictor-corrector refinement step. Overall, ItE needs less iterations than ItOTS, but the peaks are of similar magnitude. \RA{The peaks occur in the incremental steps, where the crack evolves, with around $40$ iterations on the linearizations. In incremental steps without crack evolvement, around $2$ iterations on the linearization suffice to fall below $10^{-1}$.}

Specifically on finer meshes, 
both ItL approaches yield a faster growing crack, that is closer to the true 
physics of the governing model (see again \cite{Wi17_SISC,Wi17_CMAME}), 
but, of course, at the cost of performance. Since the same incremental step is solved several times due to the ItL (and also the adaptive mesh refinement), the number of Newton active set iterations becomes very high for critical incremental steps.

\RA{Lastly, we again want to show the parallel performance by comparing the CPU time (the amount of time at least $1$ processor works) of one incremental step with no ItL on $1$ core and on $16$ cores. On $1$ core, the CPU time is $647\si{s}$ (approx. $10.6$ minutes) and on $16$ cores, it is $39\si{s}$.}
\begin{figure}[htbp!]
\begin{minipage}[b]{0.49\textwidth}
    \centering
    \includegraphics[scale = 0.26]{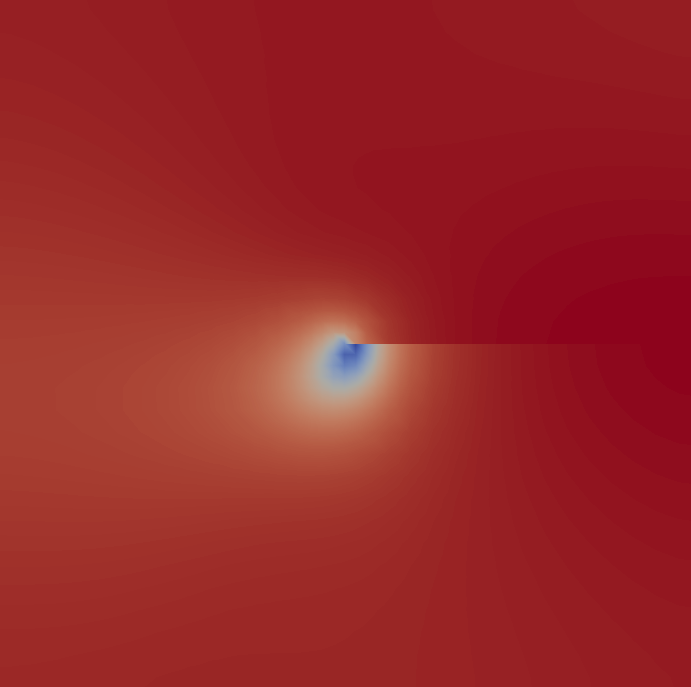}
    \caption{Visualization of the fracture for the SENS test with $h=0.022$ after $100$ incremental steps with no ItL. \RA{Red represents the fully intact part of the domain, blue the fully broken part and white stands for the transition zone.}}
    \label{fig:sens_vis_no_extra}
\end{minipage}
\hfill
\begin{minipage}[b]{0.49\textwidth}
    \centering
    \includegraphics[scale=0.26]{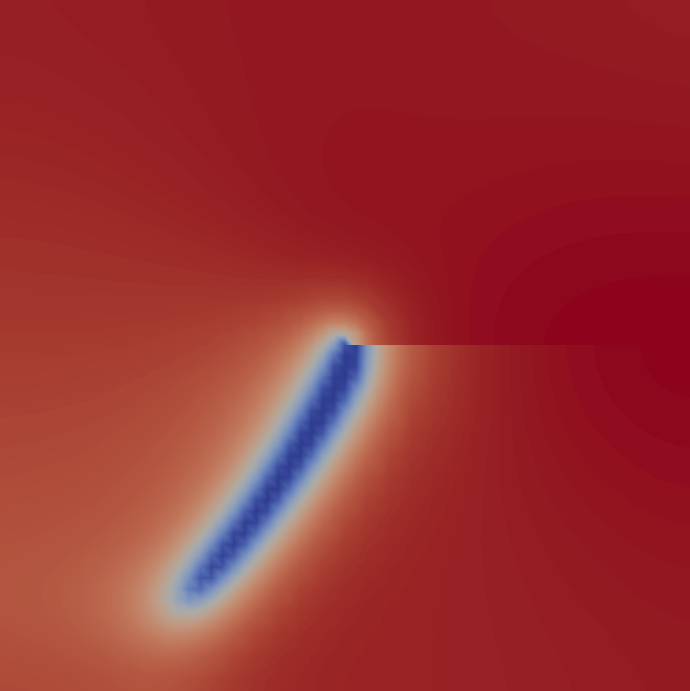}
    \caption{Visualization of the fracture for the SENS test with $h=0.022$ after $100$ incremental steps with ItE. \RA{Red represents the fully intact part of the domain, blue the fully broken part and white stands for the transition zone.}}
    \label{fig:sens_vis_extra}
\end{minipage}
\end{figure}

\begin{figure}[htbp!]
\begin{minipage}[b]{0.49\textwidth}
\scriptsize
\centering
\begin{tikzpicture}[scale = 1.0]
\begin{axis}[
    xlabel = time \auth{$[\si{s}]$}, 
    ymode = log,
    ylabel = $\#$ Active set iterations,
    legend style = {at={(0.75,1.3)}},
    mark size=0.3pt,
    grid =major,
    y post scale = 0.9,
    ]
        \addplot[line width = 6pt,color=blue]
        table[col sep=space] {miehe_shear_as_iter_case_2_0.022.txt};
        \addlegendentry{No ItL $|$ $h=0.022$}
        \addplot[line width = 6pt,color=orange,dashed]
        table[col sep=space] {miehe_shear_as_iter_extra_case_2_0.022.txt};
        \addlegendentry{ItE $|$ $h=0.022$}
        \addplot[line width = 6pt,color=violet,dotted]
        table[col sep=space]{miehe_shear_as_iter_lin_old_iter_1e-1_case_2_0.022.txt};
        \addlegendentry{ItOTS $|$ $h=0.022$}
\end{axis}
\end{tikzpicture}
\caption{Visualization of the number of active set iterations ($y$-axis) with and without iteration on the linearization (ItE and ItOTS) and $h=0.022$ depending on the time ($x$-axis) for the SENS test.}
\label{plot_as_iter_extra_iter_0.022_ms}
\end{minipage}
\hfill
\begin{minipage}[b]{0.49\textwidth}
\scriptsize
\centering
\begin{tikzpicture}[scale = 1.0]
\begin{axis}[
    xlabel = time \auth{$[\si{s}]$}, 
    ylabel = Load $F_x$ \auth{$[\si{\newton}]$},
    legend style = {at={(0.75,1.3)}},
    mark size=0.3pt,
    grid =major,
    y post scale = 0.9,
    ]
        \addplot[line width = 6pt,color=blue]
        table[col sep=space] {miehe_shear_load_disp_case_2_0.022.txt};
        \addlegendentry{No ItL $|$ $h=0.022$}
        \addplot[line width = 6pt,color=orange,dashed]
        table[col sep=space] {miehe_shear_load_disp_extra_iter_case_2_0.022.txt};
        \addlegendentry{ItE $|$ $h=0.022$}
        \addplot[line width = 6pt,color=violet,dotted]
        table[col sep=space]{miehe_shear_load_disp_lin_old_iter_1e-1_case_2_0.022.txt};
        \addlegendentry{ItOTS $|$ $h=0.022$}
\end{axis}
\end{tikzpicture}
\caption{Visualization of the $x$-load displacement (y-axis) with and without iteration on the linearization (ItE and ItOTS) and $h = 0.022$ depending on the time ($x$-axis) for the SENS test.}
\label{plot_load_disp_extra_iter_0.022_ms}
\end{minipage}
\end{figure}

\begin{figure}[htbp!]
\begin{minipage}[b]{0.49\textwidth}
\scriptsize
\centering
\begin{tikzpicture}[scale = 1.0]
\begin{axis}[
    xlabel = time \auth{$[\si{s}]$}, 
    ymode = log,
    ylabel = $\#$ Active set iterations,
    legend style = {at={(0.75,1.3)}},
    mark size=0.3pt,
    grid =major,
    y post scale = 0.9,
    ]
        \addplot[line width = 6pt,color=blue]
        table[col sep=space] {miehe_shear_as_iter_case_2_0.011.txt};
        \addlegendentry{No ItL $|$ $h=0.011$}
        \addplot[line width = 6pt,color=orange,dashed]
        table[col sep=space] {miehe_shear_as_iter_extra_case_2_0.011.txt};
        \addlegendentry{ItE $|$ $h=0.011$}
        \addplot[line width = 6pt,color=violet,dotted]
        table[col sep=space]{miehe_shear_as_iter_lin_old_iter_1e-1_case_2_0.022.txt};
        \addlegendentry{ItOTS $|$ $h=0.011$}
\end{axis}
\end{tikzpicture}
\caption{Visualization of the number of active set iterations ($y$-axis) with and without iteration on the linearization (ItE and ItOTS) and $h=0.011$ depending on the time ($x$-axis) for the SENS test.}
\label{plot_as_iter_extra_iter_0.011_ms}
\end{minipage}
\hfill
\begin{minipage}[b]{0.49\textwidth}
\scriptsize
\centering
\begin{tikzpicture}[scale = 1.0]
\begin{axis}[
    xlabel = time \auth{$[\si{s}]$}, 
    ylabel = Load $F_x$ \auth{$[\si{\newton}]$},
    legend style = {at={(0.75,1.3)}},
    mark size=0.3pt,
    grid =major,
    y post scale = 0.9,
    ]
        \addplot[line width = 6pt,color=blue]
        table[col sep=space] {miehe_shear_load_disp_case_2_0.011.txt};
        \addlegendentry{No ItL $|$ $h=0.011$}
        \addplot[line width = 6pt,color=orange,dashed]
        table[col sep=space] {miehe_shear_load_disp_extra_iter_case_2_0.011.txt};
        \addlegendentry{ItE $|$ $h=0.011$}
        \addplot[line width = 6pt,color=violet,dotted]
        table[col sep=space]{miehe_shear_load_disp_lin_old_iter_1e-1_case_2_0.011.txt};
        \addlegendentry{ItOTS $|$ $h=0.011$}
\end{axis}
\end{tikzpicture}
\caption{Visualization of the $x$-load displacement (y-axis) with and without iteration on the linearization (ItE and ItOTS) and $h = 0.011$ depending on the time ($x$-axis) for the SENS test.}
\label{plot_load_disp_extra_iter_0.011_ms}
\end{minipage}
\end{figure}

\begin{figure}[htbp!]
\begin{minipage}[b]{0.49\textwidth}
\scriptsize
\centering
\begin{tikzpicture}[scale = 1.0]
\begin{axis}[
    xlabel = time \auth{$[\si{s}]$}, 
    ymode = log,
    ylabel = $\#$ Active set iterations,
    legend style = {at={(0.75,1.3)}},
    mark size=0.3pt,
    grid =major,
    y post scale = 0.9,
    ]
        \addplot[line width = 6pt,color=blue]
        table[col sep=space] {miehe_shear_as_iter_case_2_0.0055.txt};
        \addlegendentry{No ItL $|$ $h=0.0055$}
        \addplot[line width = 6pt,color=orange,dashed]
        table[col sep=space] {miehe_shear_as_iter_extra_case_2_0.0055.txt};
        \addlegendentry{ItE $|$ $h=0.0055$}
        \addplot[line width = 6pt,color=violet,dotted]
        table[col sep=space]{miehe_shear_as_iter_lin_old_iter_1e-1_case_2_0.0055.txt};
        \addlegendentry{ItOTS $|$ $h=0.0055$}
\end{axis}
\end{tikzpicture}
\caption{Visualization of the number of active set iterations ($y$-axis) with and without iteration on the linearization (ItE and ItOTS) and $h=0.0055$ depending on the time ($x$-axis) for the SENS test.}
\label{plot_load_as_iter_extra_iter_0.0055_ms}
\end{minipage}
\hfill
\begin{minipage}[b]{0.49\textwidth}
\scriptsize
\centering
\begin{tikzpicture}[scale = 1.0]
\begin{axis}[
    xlabel = time \auth{$[\si{s}]$}, 
    ylabel = Load $F_x$ \auth{$[\si{\newton}]$},
    legend style = {at={(0.75,1.3)}},
    mark size=0.3pt,
    grid =major,
    y post scale = 0.9,
    ]
        \addplot[line width = 6pt,color=blue]
        table[col sep=space] {miehe_shear_load_disp_case_2_0.0055.txt};
        \addlegendentry{No ItL $|$ $h=0.0055$}
        \addplot[line width = 6pt,color=orange,dashed]
        table[col sep=space] {miehe_shear_load_disp_extra_iter_case_2_0.0055.txt};
        \addlegendentry{ItE $|$ $h=0.0055$}
        \addplot[line width = 6pt,color=violet,dotted]
        table[col sep=space]{miehe_shear_load_disp_lin_old_iter_1e-1_case_2_0.0055.txt};
        \addlegendentry{ItOTS $|$ $h=0.0055$}
\end{axis}
\end{tikzpicture}
\caption{Visualization of the $x$-load displacement (y-axis) with and without iteration on the linearization (ItE and ItOTS) and $h = 0.0055$ depending on the time ($x$-axis) for the SENS test.}
\label{plot_load_disp_extra_iter_0.0055_ms}
\end{minipage}
\end{figure}

\begin{figure}[htbp!]
\begin{minipage}[b]{0.49\textwidth}
\scriptsize
\centering
\begin{tikzpicture}[scale = 1.0]
\begin{axis}[
    xlabel = time \auth{$[\si{s}]$}, 
    ymode = log,
    ylabel = $\#$ Active set iterations,
    legend style = {at={(0.75,1.3)}},
    mark size=0.3pt,
    grid =major,
    y post scale = 0.9,
    ]
        \addplot[line width = 6pt,color=blue]
        table[col sep=space] {miehe_shear_as_iter_case_2_0.0027.txt};
        \addlegendentry{No ItL $|$ $h=0.0027$}
        \addplot[line width = 6pt,color=orange,dashed]
        table[col sep=space] {miehe_shear_as_iter_extra_case_2_0.0027.txt};
        \addlegendentry{ItE $|$ $h=0.0027$}
        \addplot[line width = 6pt,color=violet,dotted]
        table[col sep=space]{miehe_shear_as_iter_lin_old_iter_1e-1_case_2_0.0027.txt};
        \addlegendentry{ItOTS $|$ $h=0.0027$}
\end{axis}
\end{tikzpicture}
\caption{Visualization of the number of active set iterations ($y$-axis) with and without iteration on the linearization (ItE and ItOTS) and $h=0.0027$ depending on the time ($x$-axis) for the SENS test.}
\label{plot_load_as_iter_extra_iter_0.0027_ms}
\end{minipage}
\hfill
\begin{minipage}[b]{0.49\textwidth}
\scriptsize
\centering
\begin{tikzpicture}[scale = 1.0]
\begin{axis}[
    xlabel = time \auth{$[\si{s}]$}, 
    xmax = 0.0236,
    ylabel = Load $F_x$ \auth{$[\si{\newton}]$},
    legend style = {at={(0.75,1.3)}},
    mark size=0.3pt,
    grid =major,
    y post scale = 0.9,
    ]
        \addplot[line width = 6pt,color=blue]
        table[col sep=space] {miehe_shear_load_disp_case_2_0.0027.txt};
        \addlegendentry{No ItL $|$ $h=0.0027$}
        \addplot[line width = 6pt,color=orange,dashed]
        table[col sep=space] {miehe_shear_load_disp_extra_iter_case_2_0.0027.txt};
        \addlegendentry{ItE $|$ $h=0.0027$}
        \addplot[line width = 6pt,color=violet,dotted]
        table[col sep=space]{miehe_shear_load_disp_lin_old_iter_1e-1_case_2_0.0027.txt};
        \addlegendentry{ItOTS $|$ $h=0.0027$}
\end{axis}
\end{tikzpicture}
\caption{Visualization of the $x$-load displacement (y-axis) with and without iteration on the linearization (ItE and ItOTS) and $h = 0.0027$ depending on the time ($x$-axis) for the SENS test.}
\label{plot_load_disp_extra_iter_0.0027_ms}
\end{minipage}
\end{figure}

\newpage
\section{Conclusion}\label{sec_conclusions}
In this work, we investigated and improved the primal-dual active set phase-field fracture 
formulation 
proposed in \cite{heister2015primal}. From the three major advancements is one of theoretical 
nature and two are numerically motivated. 
On the theoretical side we worked out details
of deriving the active set formulation from the governing complementarity system.
In numerics, we examined the active set constant $c$ and proposed four different 
test cases to be compared. Thirdly, we studied an iteration on the linearization 
procedure of the phase-field variable in the displacement equation.
Concerning the derivations established in Section \ref{sec:problem_formulation}, we 
started from the coupled variational inequality system (CVIS) and derived the 
complementarity formulation in weak form. Then, this was linked to the corresponding strong form 
conditions, which are the starting point for the active set algorithm.
In Section \ref{sec:numerical_solution}, we investigated theoretically the role of the active set constant $c$ in the solution algorithm and found a new setting, which avoids the wrong classification of degrees of freedom. This increases the convergence speed of the Newton active set method on fine meshes and thus yields better performances. 
In Section \ref{sec:iteration_on_extra}, besides studying the iteration on the linearization,
we formulated a final new algorithm.
For both linearization approaches, ItL resolves the time-lagging issue, 
that is due to the resulting discretization error.
\auth{In Section \ref{sec:numerical_tests}, we performed five numerical tests 
including two- and three-dimensional settings, as well as stationary crack that only 
vary in their width, and propagating fractures.
In the SENS test (Section \ref{subsec:sens})}, ItE performs better than ItOTS in terms of fracture evolution speed, but ItE shows some non-physical behaviour due to the iterations 
on the extrapolation. In terms of accuracy, both approaches 
yield findings that are closer to the expected physics, as it can be compared 
with our prior work when using a fully monolithic model.
Our overall conclusion of both modifications, namely active set constant and 
the ItL schemes, is that one has to make a choice between accuracy 
and efficiency. This is not surprising at all, but highlighted in this work for our primal-dual phase-field fracture framework.

\section*{Acknowledgements}
All authors thank Viktor Kosin (Universit\'e Paris-Saclay) for fruitful discussions 
on the iteration on the extrapolation part and Johannes Lankeit (Leibniz University Hannover) for giving valuable advice and comments on solution spaces in connection with equivalence proofs within this work. Moreover, the authors thank Sebastian Bohlmann for the Scientific Computing environment at IfAM. The present work has been \auth{partially} carried out within the DFG Collaborative Research Center (CRC) 1463 “Integrated design and operation methodology for offshore megastructures”, which is funded by the Deutsche Forschungsgemeinschaft (DFG, German Research Foundation) - Project-ID 434502799, SFB 1463. \auth{Moreover, the authors thank the (anonymous) reviewers for their questions that helped to improve the manuscript.}


\end{document}